%% file: Infty-Bundles_main.tex
\documentclass[twoside, 11pt]{article}

\usepackage[latin1]{inputenc}
\usepackage[T1]{fontenc}
\usepackage[english]{babel}
\usepackage{color}
\usepackage{tocloft}
\usepackage{verbatim}

\usepackage{tikz}
\usetikzlibrary{snakes}
\usetikzlibrary{arrows}
\usepackage{tikz-cd}

\usepackage[show]{optional}

\usepackage[
  hmarginratio={1:1},     
  vmarginratio={1:1},     
  textwidth=17cm,        
  textheight=24cm,
  heightrounded          
]{geometry}

\usepackage[noindentafter]{titlesec}
\titlespacing{\paragraph}{%
  0pt}{
0.25\baselineskip}{
1em}

\titlespacing{\section}{%
0pt}{
0.2cm}{
0em}

\titlespacing{\subsection}{%
0pt}{
0.2cm}{
0em}

\titlespacing{\subsubsection}{%
0pt}{
0cm}{
0em}

\usepackage{amsmath,amsthm,amssymb}
\usepackage{mathrsfs, mathtools, xpatch}
\usepackage{dsfont}

\allowdisplaybreaks

\usepackage[all]{xy} 
\SelectTips{cm}{11} 
\usepackage[hidelinks]{hyperref}

\usepackage[numbers,sort&compress]{natbib}
\makeatletter
\renewcommand{\@biblabel}[1]{[#1]\hfill}
\makeatother
\usepackage[mathscr]{euscript}
\usepackage{calrsfs}
\DeclareMathAlphabet{\pazocal}{OMS}{zplm}{m}{n}

\usepackage[shortlabels]{enumitem}

\newtheoremstyle{thm}                                                           
{0.15cm}                                         
{0.15cm}                                         
{\itshape}      
{}                                      
{\bfseries}                             
{}                                      
{0.2cm}                                         
{\thmname{#1}~\thmnumber{#2}\thmnote{ (#3)}}%

\makeatletter
\xpatchcmd{\proof}{\topsep6\p@\@plus6\p@\relax}{}{}{}
\makeatother

\newtheoremstyle{rmk}                                                           
{0.15cm}                                         
{0.15cm}                                         
{}      
{}                                      
{\bfseries}                             
{}                                      
{0.2cm}                                         
{}                                      

\theoremstyle{thm}
\newtheorem{theorem}[equation]{Theorem}

\newtheorem{corollary}[equation]{Corollary}

\newtheorem{proposition}[equation]{Proposition}

\theoremstyle{rmk}
\newtheorem{example}[equation]{Example}
\newtheorem{remark}[equation]{Remark}
\newtheorem{definition}[equation]{Definition}

\numberwithin{equation}{section}

\DeclareMathAlphabet{\mathbbe}{U}{bbold}{m}{n}

\input{Infty-Bundles_Defs.tex}

\mathtoolsset{showonlyrefs,showmanualtags}

\newlength{\Displayskip}
\begin{document}

\setlength{\abovedisplayskip}{\Displayskip}
\setlength{\belowdisplayskip}{\Displayskip}

\begin{center}
\LARGE{\textbf{
$\infty$-Bundles
}}
\end{center}
\vspace{0.2cm}
\begin{center}
\large Severin Bunk
\end{center}

\vspace{0.2cm}

\begin{abstract}
\noindent
Higher bundles are homotopy coherent generalisations of classical fibre bundles.
They appear in numerous contexts in geometry, topology and physics.
In particular, higher principal bundles provide the geometric framework for higher-group gauge theories with higher-form gauge potentials and their higher-dimensional holonomies.
An $\infty$-categorical formulation of higher bundles further allows one to identify these objects in contexts outside the worlds of smooth manifolds or topological spaces.
This article reviews the theory of $\infty$-bundles, focussing on principal $\infty$-bundles, and surveys several of their applications.
It is an invited contribution to the Topology section in the second edition of the Encyclopedia of Mathematical Physics.
\end{abstract}



\tableofcontents


\section{Introduction}
\label{sec:Introduction}


One of the most influential and successful paradigms in mathematics in recent decades has been to weaken and generalise familiar objects by replacing structural identities with homotopy coherence data.
This allows one to recognise (these new versions of) familiar structures in new, previously inaccessible places.
In this homotopy coherent sense, the based loop space of a pointed topological space, together with path concatenation, is in fact a group, and the cochain complex of a topological space, together with the cup product and the Eilenberg-Zilber map, is a commutative ring.

The main subject of this article is the theory of $\infty$-bundles, the corresponding homotopy coherent weakening of the study of fibre bundles.
We emphasise the notion of principal $\infty$-bundles, as this is in some ways the most fundamental type of higher bundles.
Here, essentially, two classical concepts need to be replaced by homotopically richer counterparts:
the concept of a surjective map from the total space of a principal bundle to the base space needs to be replaced by a map which is `homotopically surjective' (concretely: an effective epimorphism), and the concept of a group needs to be replaced by a higher categorical object with a homotopy coherent multiplication in which everything is invertible (concretely: a group object in an $\infty$-category).
Historically, achieving these replacements in a fully general setting relied crucially on developments in the theory of $\infty$-categories~\cite{BV:Ho-invar_alg_strs, Joyal:QCats_and_applications, Lurie:HTT}.

In geometry and mathematical physics principal bundles underpin the study of gauge fields, or connections.
Gauge theory and its geometry have long been protagonists in the lasting synergy between mathematics and physics (see, for instance,~\cite{Nash:DiffTop_and_QFT} for a range of results in this direction).
Early mathematical appearances%
\footnote{This introduction is not intended to provide a complete guide to the literature, and references will necessarily be incomplete.
The author apologises for any omissions.}
of higher bundles include~\cite{Giraud:Cohomologie_non-abelienne, Breen:Classif_of_2-gerbes_and_2-stacks, Brylinski:Loops_and_GeoQuan}, seeking for geometric descriptions of non-abelian and higher differential cohomology.
The case of integer differential cohomology was strongly advanced by Gajer in~\cite{Gajer:Geo_of_Deligne_Coho}, who established a first notion of higher $\rmU(1)$-bundles and connections as cocycles for differential cohomology in any degree.
In physics interest in higher bundles was sparked by string and M-theory, which naturally posed the question for two-dimensional generalisations of parallel transport and, thus, 2-form connections on principal bundles.
This two-dimensional case was first formalised in the influential papers~\cite{BS:HGT_2-conns_on_2-bdls, BS:HGT}.
Here one considers categorified principal bundles whose structure group is a 2-group.
In the most general case, this is a monoidal groupoid whose objects are each invertible with respect to the monoidal structure, but often stricter models, such as crossed modules were used~\cite{ACJ:Nonabelian_BGrbs, NW:Nonabelian_gerbes}.
Subsequently, these concepts were developed further to classify principal bundles for structure groups with increasingly higher structure and thus richer homotopy coherence data (see, for instance,~\cite{JL:Higher_Principal_Bdls, RS:Spl_pr_bdls_in_parameterised_spaces, JSW:Semistrict_HGT}).
The full definition of principal $\infty$-bundles, with the weakened concepts of surjections and groups mentioned above, first appeared in~\cite{NSS:oo-bundles_I}.

This paper reviews that notion of principal $\infty$-bundles and some of its applications in geometry, topology and mathematical physics.
We cover the essential background, key definitions and results, and give examples throughout, which we hope will help convey the intuition behind the theory as well as signpost the reader to further applications.
In Section~\ref{sec:Motivation} we motivate the passage from ordinary to $\infty$-bundles with an illustrating example.
For the reader's convenience, we include a very brief sketch of $\infty$-categories in Section~\ref{sec:oo-cats}.
In Sections~\ref{sec:EEpis} and~\ref{sec:Groups and actions} we review the particular class of $\infty$-categories in which a theory of $\infty$-bundles can be formulated---the so-called $\infty$-topoi---and recall the notions of groups and group actions in $\infty$-topoi.
Section~\ref{sec:Pr-oo-bundles} covers the definition and characterisations of principal $\infty$-bundles in the sense of~\cite{NSS:oo-bundles_I}.
An important feature of this theory is the existence of classifying objects for principal $\infty$-bundles, which we review in Section~\ref{sec:Non-ab coho and BG}.
In Section~\ref{sec:associated bundles} we recall how essentially every $\infty$-bundle can be obtained as an associated bundle for some principal $\infty$-bundle, and in Section~\ref{sec:Geometric mps} we recall that principal $\infty$-bundles are preserved under $\infty$-functors with certain properties.
Section~\ref{sec:Connections} contains a very brief survey of the still incomplete theory of connections on $\infty$-bundles.
We close in Section~\ref{sec:Appls in physics} by outlining three applications of higher bundles in mathematical physics.


\paragraph*{Acknowledgements}


The author would like to thank L.~Müller, J.~Nuiten, C.~Sämann, U.~Schreiber, R.~Szabo and K.~Waldorf for many enlightening discussions about higher bundles.
The author is supported by the Deutsche Forschungsgemeinschaft (DFG, German Research Foundation) under project number 468806966.


\section{Towards homotopy coherent principal bundles}
\label{sec:Motivation}


Let $\Top$ denote the category of topological spaces and continuous maps.
Recall that, in $\Top$, principal bundles are defined as follows:

\begin{definition}
\label{def:ord principal bundle}
Let $G$ be a topological group.
A principal $G$-bundle on a topological space $X \in \Top$ consists of a continuous map $p \colon P \to X$ and a continuous right $G$-action on $P$ such that
\begin{myenumerate}
\item (\textit{Fibre-preserving action}) the $G$-action preserves the fibres of the map $p \colon P \to X$,

\item (\textit{Local triviality}) there exists an open covering $\CU = \{U_a\}_{a \in \Lambda}$ of $X$ and homeomorphisms $\varphi_a \colon P_{|U_a} \to U_a \times G$ which intertwine the right action of $G$ on $P$ with the canonical action of $G$ on itself via right multiplication and which commute with the projections to $U_a$.

\item (\textit{Principality condition})
the \textit{shear map} $P \times G \to P \times_X P$, $(p,g) \mapsto (p, pg)$ is a homeomorphism.
\end{myenumerate}
A \textit{morphism of principal $G$-bundles $(P \to X) \longrightarrow (Q \to X)$ over $X$} is a fibre-preserving morphism $f \in \Top_{/X}(P,Q)$ which commutes with the $G$-actions.
\end{definition}

\begin{remark}
\label{rmk:classical pr-Bun mps are isos}
One can show that each morphism of principal $G$-bundles on $X \in \Top$ is automatically an \textit{isomorphism}.
\qen
\end{remark}

For each principal $G$-bundle $p \colon P \to X$, there is a canonical homeomorphism $P/G \cong X$.
Because of the principality condition, we may also write
\begin{equation}
	P/G \cong \colim
	\begin{tikzcd}
		\big( P \underset{X}{\times} P \ar[r, shift left=0.1cm] \ar[r, shift left=-0.1cm]
		& P \big)\,. \ar[l]
	\end{tikzcd}
\end{equation}

However, there exist maps $p\colon P \to X$ in $\Top$ which one would like to consider as instances of principal bundles, but which do not fit into Definition~\ref{def:ord principal bundle}.
For instance, let $X \in \Top$ be connected and fix a basepoint $x_0 \in X$.
We let $P_{x_0} X$ be the based path space of $(X, x_0)$ and $\ev_1 \colon P_{x_0} X \to X$ the endpoint evaluation map, sending a path $\gamma \colon [0,1] \to X$ to its endpoint $\gamma(1)$.
There is a canonical equivalence
\begin{equation}
\label{eq:shear map for based path fibration}
	P_{x_0} X \underset{X}{\times} P_{x_0} X \longrightarrow \Omega_{x_0} X\,,
	\qquad
	(\gamma_0, \gamma_1) \longmapsto \overline{\gamma_1} * \gamma_0\,,
\end{equation}
where $\Omega_{x_0} X$ is the based loop space of $(X, x_0)$, $(-)*(-)$ denotes path concatenation, and $\overline{(-)}$ denotes path reversal.
It is a classical result that the based loop space $\Omega_{x_0} X$ is not simply a space, but carries a homotopy coherent weakening of a group structure.
It is a grouplike $\bbE_1$-algebra, or $A_\infty$-algebra, in $\Top$~\cite{Stasheff:HoAss_of_H-spaces_I, Stasheff:HoAss_of_H-spaces_II, May:Geo_of_iterated_loop_spaces}.
We also have a basepoint preserving, homotopy coherent action map
\begin{equation}
	P_{x_0} X \times \Omega_{x_0} X \longrightarrow P_{x_0} X\,,
	\qquad
	(\gamma, \alpha) \longmapsto \gamma * \alpha\,,
\end{equation}
which makes the shear map
\begin{equation}
	P_{x_0} X \times \Omega_{x_0} X \longrightarrow P_{x_0} X \underset{X}{\times} P_{x_0} X\,,
	\qquad
	(\gamma, \alpha) \longmapsto \big( \gamma, \gamma * \alpha)
\end{equation}
not into a homeomorphism, but into a homotopy equivalence.

Moreover, the map $\ev_1 \colon P_{x_0} X \to X$ is locally trivial in the following sense:
let $\CU = \{U_a\}_{a \in \Lambda}$ be an open covering of $X$ such that each $U_a$ and each finite intersection $U_{a_0 \cdots a_n} \coloneqq U_{a_0} \cap \cdots \cap U_{a_n}$ is a contractible open subset of $X$.
The contractibility of $U_a$ means, in particular, that we can find a section $\Gamma_a \colon U_a \to (P_{x_0} X)_{|U_a}$ of $\ev_1$ over $U_a \subset X$, for each $a \in \Lambda$.
We thus obtain (again, not a homeomorphism, but) a homotopy equivalence
\begin{equation}
	\varphi_a \colon (P_{x_0} X)_{|U_a} \longrightarrow U_a \times \Omega_{x_0} X\,,
	\qquad
	\gamma \longmapsto \big( \gamma(1), \overline{\Gamma_a(\gamma(1))} * \gamma \big)\,.
\end{equation}
Finally, from these trivialisations we obtain transition maps
\begin{equation}
	\varphi_b \circ \varphi_a^{-1}
	= \big( \id_{U_{ab}}, (-) * \Gamma_{ab} \big)
	\colon U_{ab}\times \Omega_{x_0} X \longrightarrow U_{ab} \times \Omega_{x_0} X\,,
\end{equation}
which indeed consist of the identity on $U_{ab}$ together with the right action of the transition functions
\begin{equation}
	\Gamma_{ab} \colon U_{ab} \longrightarrow \Omega_{x_0} X\,,
	\qquad
	\Gamma_{ab}(x) = \overline{\Gamma_b(x)} * \Gamma_a(x)\,,
\end{equation}
as we would expect from transition maps for principal bundles.
However, instead of a cocycle \textit{identity} for the group-valued functions $\Gamma_{ab}$, we only have a canonical homotopy, i.e.~a 1-simplex in the space of maps $U_{abc} \to \Omega_{x_0} X$,
\begin{equation}
	 \Gamma_{abc} \colon \Gamma_{bc} * \Gamma_{ab} \eq \Gamma_{ac}\,,
\end{equation}
for each $a,b,c \in \Lambda$, and further homotopies between homotopies over quadruple overlaps $U_{abcd}$, and so on.
It is a crucial step in passing from ordinary to higher bundles that one needs to view all these homotopies as part of the transition data of the bundle $P_{x_0} X \to X$.

This perspective on the path fibration $\ev_1 \colon P_{x_0} X \to X$ motivates the search for an abstraction of the concept of a principal bundle which is sufficiently weakened, or categorified, to capture this example.
Ideally, we need to formulate any such generalised version of principal bundles in a fashion which does not make explicit reference to the category $\Top$, but rather extracts only those of the many properties of the category $\Top$ which we actually need.
That will allow us to consider principal bundles in much more general contexts, such as smooth spaces with higher structure and possibly derived structure, or even in completely different situations.
At the very least, in a context where principal bundles in this sense would be feasible, we need to have notions of
\begin{itemize}[itemsep=-0.1cm, leftmargin=*, topsep=0cm]
\item covering, or surjection,

\item weakly and coherently invertible morphisms, or equivalences,

\item homotopy coherent groups and group actions,

\item quotients by group actions, and

\item pullbacks along any morphism.
\end{itemize}
The biggest conceptual leap consists in encoding homotopy coherence.
This is most conveniently achieved by passing to the framework of $\infty$-categories.
It turns out that there is a class of $\infty$-categories that provides a particularly well-adapted background for this enhanced theory of principal bundles.
It consists of the \textit{$\infty$-topoi}, whose definition we build up to in the next couple of sections (Definition~\ref{def:oo-topos}), before giving the $\infty$-categorical definition of principal bundles in Definition~\ref{def:G-pr oo-Bun}.


\section{From categories to $\infty$-categories}
\label{sec:oo-cats}


We give a short---and by no means complete---introduction to $\infty$-categories.
By the term $\infty$-category we shall always mean an $(\infty,1)$-category, i.e.~we allow for non-trivial $k$-morphisms for each $k \in \NN$, but the morphisms in levels $k > 1$ are invertible.
Explicitly, we model $\infty$-categories as simplicial sets satisfying the inner horn-lifting conditions (also known as \textit{quasicategories}).
We refer the reader to~\cite{Joyal:QCats_and_applications, Lurie:HTT, Cisinski:HCats_HoAlg} for comprehensive treatments of $\infty$-category theory in this language.

Let $\bbDelta$ denote the category of finite, totally ordered sets of the form $[k] = \{0, 1, \ldots, k\}$, for $k \in \NN_0$, and order-preserving maps.
The category of \textit{simplicial sets} is the functor category
\begin{equation}
	\sSet \coloneqq \Fun(\bbDelta^\opp, \Set)\,.
\end{equation}
For each $n \in \NN_0$, there is a standard $n$-simplex $\Delta^n = \bbDelta(-, [n]) \in \sSet$, and for each $0 \leq k \leq n$ there is the $k$-th horn $\Lambda^n_k \subset \Delta^n$ obtained, pictorially, by removing the interior of the $n$-simplex and the $(n{-}1)$-dimensional face opposite the $k$-th vertex.

\begin{definition}
\cite{BV:Ho-invar_alg_strs}
An \textit{$\infty$-category}, or \textit{quasicategory}, is a simplicial set $\scC \in \sSet$ satisfying that each diagram of solid arrows
\begin{equation}
\begin{tikzcd}
	\Lambda^n_k \ar[r] \ar[d, hookrightarrow]
	& \scC
	\\
	\Delta^n \ar[ur, dashed]
	&
\end{tikzcd}
\end{equation}
admits a lift as indicated, for each $n \geq 2$ and $0 < k < n$.
A \textit{morphism of $\infty$-categories $\scC \to \scD$}, or \textit{$\infty$-functor}, is a morphism of simplicial sets $\scC \to \scD$.
\end{definition}

\begin{example}
Given a 1-category $C$, its \textit{nerve} $NC \in \sSet$ is defined by setting
\begin{equation}
	NC_n = \Fun([n],C)\,,
	\qquad
	\forall\ n \in \NN_0\,.
\end{equation}
This is always an $\infty$-category, and we obtain a fully faithful inclusion of the (2,1)-category of 1-categories, functors and natural isomorphisms into the $\infty$-category of $\infty$-categories (defined in Example~\ref{eg:localisations}(2) below).
We have that $N[n] = \Delta^n$, for each $n \in \NN_0$.
\qen
\end{example}

\begin{example}
Given an $\infty$-category $\scC$ and any simplicial set $K \in \sSet$, the internal hom $\scC^K \in \sSet$ (with $(\scC^K)_n = \sSet(K \times \Delta^n, \scC)$) is again an $\infty$-category.
For $\scD$ another $\infty$-category, the \textit{$\infty$-category of functors $\scC \to \scD$} is the internal hom $\scFun(\scC, \scD) \coloneqq \scD^\scC$ in $\sSet$.
\qen
\end{example}

In contrast to ordinary (1-)categories, in an $\infty$-category there is no composition law.
Instead, we can interpret the $n$-simplices $\sigma \in \scC_n = \sSet(\Delta^n, \scC)$ in an $\infty$-category $\scC$ as encoding how a given 1-morphism (the image of the edge $0 \to n$ in $\Delta^n$) can be written as compositions of $n$ many 1-morphisms, which are given as the image of the \textit{spine}
\begin{equation}
\begin{tikzcd}
	\Sp^n = \Delta^{\{0,1\}} \underset{\Delta^{\{1\}}}{\sqcup} \cdots \underset{\Delta^{\{n-1\}}}{\sqcup} \Delta^{\{n-1,n\}}
	\ar[r, hookrightarrow]
	& \Delta^n\,.
\end{tikzcd}
\end{equation}

For each such sequence of $n$ composable 1-morphisms, there is a contractible space of choices for such composition data:

\begin{proposition}
\emph{\cite[Cor.~3.7.6]{Cisinski:HCats_HoAlg}}
Given an $\infty$-category $\scC$ and $n \geq 2$, the inclusion $\Sp^n \hookrightarrow \Delta^n$ induces a trivial Kan fibration of simplicially enriched homs,
\begin{equation}
	\scFun(N[n], \scC) = \scC^{\Delta^n}
	\longrightarrow \scC^{\Sp^n} = \scFun(\Sp^n, \scC)\,.
\end{equation}
\end{proposition}

The following is a manifestation of Grothendieck's homotopy hypothesis:

\begin{definition}
An \textit{$\infty$-groupoid}, or \textit{Kan complex}, is a simplicial set $\scK \in \sSet$ satisfying that each diagram in $\sSet$ of solid arrows
\begin{equation}
\begin{tikzcd}
	\Lambda^n_k \ar[r] \ar[d, hookrightarrow]
	& \scK
	\\
	\Delta^n \ar[ur, dashed]
	&
\end{tikzcd}
\end{equation}
admits a lift as indicated, for each $n \geq 1$ and $0 \leq k \leq n$.
\end{definition}

\begin{example}
Let $G \in \Fun(\bbDelta^\opp, \Grp)$ be a simplicial group.
Then, its underlying simplicial set is a Kan complex (see, for instance,~\cite[Lemma~I.3.4]{GJ:Spl_HoThy}).
\qen
\end{example}

\begin{remark}
The nerve $NC$ of a category $C$ is a Kan complex if and only if $C$ is a groupoid.
\qen
\end{remark}

Invertible 1-morphisms in an $\infty$-category are called \textit{equivalences}.
One can show that an $\infty$-category is an $\infty$-groupoid if and only if all its morphisms are equivalences~\cite[Thm.~3.5.1]{Cisinski:HCats_HoAlg}.
All 1-categorical concepts, including slice categories, (co)limits, Kan extensions, adjunctions, equivalences and monoidal structures, have $\infty$-categorical generalisations.
These are often much richer than the 1-categorical versions (for instance, `uniqueness' generally translates to the existence of a contractible space of choices).
Importantly, this also applies to the concept of localisation:
in the $\infty$-categorical setting, this is a procedure for adding inverses for a chosen class of 1-morphisms, with additional 2-cells witnessing that these new morphisms are indeed inverse to the given 1-morphisms, and additional higher cells witnessing the necessary new higher coherences.

\begin{definition}
\cite[Def.~7.1.2]{Cisinski:HCats_HoAlg}
Let $\scC$ be an $\infty$-category and $W \subset \scC$ any simplicial subset.
An \textit{$\infty$-categorical localisation} of $\scC$ at $W$ consists of an $\infty$-category $L_W \scC$ and an $\infty$-functor $\scC \to L_W \scC$, satisfying the following properties:
\begin{myenumerate}
\item for each 1-simplex $f$ in $W$, the image $\gamma(f)$ is an equivalence in $L_W \scC$, and

\item for each $\infty$-category $\scD$, the morphism 
\begin{equation}
	\gamma^* \colon \scFun(L_W\scC, \scD) \longrightarrow \scFun_W(\scC, \scD)
\end{equation}
induces an equivalence between the $\infty$-categories of $\infty$-functors $L_W\scC \to \scD$ and the full $\infty$-subcategory of $\scFun(\scC, \scD)$ on those $\infty$-functors which send all morphisms in $W$ to equivalences.
\end{myenumerate}
\end{definition}

\begin{proposition}
{\emph{\cite[Prop.~7.1.3]{Cisinski:HCats_HoAlg}}}
Let $\scC$ be an $\infty$-category and $W \subset \scC$ any simplicial subset.
The $\infty$-categorical localisation $L_W\scC$ exists and is essentially unique.
\end{proposition}

\begin{example}
\label{eg:localisations}
Localisation allows us to construct many important examples of $\infty$-categories%
\footnote{It is even true that \textit{every} $\infty$-category $\scC$ is equivalent to an $\infty$-categorical localisation of the nerve of its 1-category $\bbDelta_{/\scC}$ of simplices; see, for instance,~\cite[Prop.~7.3.15]{Cisinski:HCats_HoAlg}.}.
Knowing how to write an $\infty$-category as a localisation of a 1-category is often helpful for explicit constructions and computations.
\begin{myenumerate}
\item Our central example is the following:
if $(C, W)$ is a relative category with weak equivalences $W$, we obtain an $\infty$-categorical localisation $L_W NC$ of the nerve of $C$ at $W$.
A rich supply for such pairs $(C,W)$ arises from model structures on $C$, where $W$ are the weak equivalences in the model structure.
For simplicial model structures in particular, there are well-controlled descriptions of the $\infty$-categorical localisations $L_W NC$ (see, for instance,~\cite[Prop.~A.3.7.6]{Lurie:HTT}).

\item There is a unique model structure on $\sSet$ whose cofibrations are the levelwise injections and whose fibrant objects are the $\infty$-categories:
this is the \textit{Joyal model structure}.
We denote its weak equivalences by $W_J$.
The \textit{$\infty$-category of $\infty$-categories} is the localisation
\begin{equation}
	\Cat_\infty \coloneqq L_{W_J} N \sSet\,.
\end{equation}

\item Similarly, we have the Kan-Quillen model structure on $\sSet$, whose cofibrations are the levelwise injections and whose fibrant objects are the Kan complexes.
Its weak equivalences are the weak homotopy equivalences of simplicial sets.
We denote this class of morphisms by $W_{KQ}$.
This gives rise to the \textit{$\infty$-category of spaces}, or \textit{$\infty$-groupoids},
\begin{equation}
	\scS \coloneqq L_{W_{KQ}} N \sSet\,.
\end{equation}

\item Given an $\infty$-category $\scC$, the $\infty$-functors $\scC^\opp \to \scS$ are called \textit{$\infty$-presheaves} on $\scC$.
We write
\begin{equation}
	\PSh(\scC) \coloneqq \scFun(\scC^\opp, \scS)\,.
\end{equation}
Given a category $C$, let $W$ denote the class of objectwise weak homotopy equivalences in the 1-category $\Fun(C^\opp, \sSet)$.
We have an equivalence of $\infty$-categories
\begin{equation}
\label{eq:presentation of PSh}
	\PSh(NC^\opp, \scS)
	\simeq L_W N \Fun(C^\opp, \sSet)
\end{equation}
(see, for instance,~\cite[Ch.~2]{Lurie:HTT}, or subsequent simplifications~\cite{HM:Left_fibs_and_hocolims, Bunk:oo-Loc}).

\item If $(C,\tau)$ is a category with a Grothendieck coverage, we can also localise model categories of homotopy sheaves of simplicial sets at their weak equivalences $W_\tau$.
This gives a presentation,
\begin{equation}
\label{eq:presentation of Sh}
	\Sh(NC,\tau) \simeq L_{W_\tau} \Fun(C^\opp, \sSet)
\end{equation}
for the $\infty$-category of $\infty$-sheaves on $NC$ (see~\cite[Sec.~6.2]{Lurie:HTT} for more background).
\qen
\end{myenumerate}
\end{example}

\begin{example}
\label{eg:oo-topoi of smooth spaces}
Example~\ref{eg:localisations}(4) and~(5) are particularly relevant to physics, because they allow us to describe $\infty$-categories of smooth spaces that generalise manifolds.
Let $\Cart$ be the category of \textit{cartesian spaces}:
its objects are all submanifolds $c \subset \RN^\infty$ such that there exists an $n \in \NN_0$ and a diffeomorphism $c \cong \RN^n$.
Its morphisms are all smooth maps between these manifolds.
The category $\Cart$ carries a Grothendieck coverage $\tau$, given by good open coverings, i.e.~those open coverings where all finite non-empty intersections of patches are again cartesian spaces.
We thus obtain $\infty$-categories
\begin{align}
	\scH &\coloneqq \PSh(N\Cart)
	\simeq L_W N\Fun(\Cart^\opp, \sSet)\,,
	\\
	\scH_\tau &\coloneqq \Sh(N\Cart, \tau)
	\simeq L_{W_\tau} N\Fun(\Cart^\opp, \sSet)\,.
\end{align}
These $\infty$-categories are frequently used to describe higher-geometric generalisations of manifolds and sheaves, in particular in contexts in mathematical physics~\cite{Schreiber:DCCT, FSS:Cech_diff_char_classes_via_L_infty}%
\footnote{See also~\cite{Nuiten:Stacks_and_fractions, Nuiten:Hypercovers} for the fact that $\scH_\tau$ contains the $\infty$-categories obtained by localising categories of Lie $n$-groupoids at the Morita equivalences (using that $(\Cart, \tau)$ has enough points).}.
There are also further enhancements of these $\infty$-categories of $\infty$-(pre)sheaves which detect infinitesimal, or derived, aspects of higher smooth spaces (see, for instance,~\cite{Schreiber:DCCT, Lurie:DAG_V:Struct_spaces, Nuiten:Thesis, Steffens:Derived_C^oo-geomtery_I, AY:Non-pert_BV-theory}).
\qen
\end{example}


\section{From surjections to effective epimorphisms and $\infty$-bundles}
\label{sec:EEpis}


We now ask for a homotopy coherent, i.e.~$\infty$-categorical, version of Definition~\ref{def:ord principal bundle}.
In the $\infty$-category $\scS$ of spaces it is no longer useful to think of a map $p \colon P \to X$ as `surjective';
instead, we should replace this concept with that of an \textit{essentially surjective} functor of $\infty$-groupoids.
A morphism $p \colon P \to X$ in $\scS$ satisfies this property if and only if the induced map of sets $\pi_0 p \colon \pi_0 P \to \pi_0 X$ is surjective.

An alternative way of saying that a continuous map $p \colon P \to X$ is surjective is that in its image factorisation
\begin{equation}
\begin{tikzcd}
	P \ar[rr, "p"] \ar[dr]
	& & X
	\\
	& \im(p) \ar[ur, hookrightarrow, "\iota"']
\end{tikzcd}
\end{equation}
the inclusion $\iota$ is a homeomorphism.
Importantly, we can write
\begin{equation}
	\im(p) = \colim
		\begin{tikzcd}
			\big( P \underset{X}{\times} P \ar[r, shift left=0.1cm] \ar[r, shift left=-0.1cm]
			& P \big)\,. \ar[l]
		\end{tikzcd}
\end{equation}
It is the quotient of $P$ by the relation that $y_0 \sim y_1$ in $P$ if and only if $p(y_0) = p(y_1)$.

This has an $\infty$-categorical enhancement:
let $\bbDelta_+ = \bbDelta^\triangleleft$ be the category $\bbDelta$ with an initial object $[-1]$ adjoined.
Let $\bbDelta_{+, \leq 0} \subset \bbDelta_+$ be the full subcategory on the objects $[-1]$ and $[0]$ (its only non-identity morphism is $[-1] \to [0]$).
Finally, let $\jmath \colon \bbDelta_{+, \leq 0} \to \bbDelta_+$ be the canonical inclusion functor.

\begin{definition}
\label{def:Cech nerves}
Let $\scC$ be an $\infty$-category with finite limits.
The \textit{augmented \v{C}ech nerve} $\cC_+p$ of a morphism $p \colon P \to X$ in $\scC$ is the right Kan extension
\begin{equation}
\begin{tikzcd}[column sep=1.25cm, row sep=0.75cm]
	N \bbDelta_{+, \leq 0}^\opp \ar[d, hookrightarrow, "\jmath"'] \ar[r, "\{p\}"]
	& \scC
	\\
	N\bbDelta_+^\opp \ar[ru, dashed, "\cC_+ p = \jmath_*\{p\}"']
	&
\end{tikzcd}
\end{equation}
The \textit{\v{C}ech nerve} of $p \colon P \to X$ is the simplicial object $\cC p \colon N\bbDelta^\opp \to \scC$ underlying $\cC_+ p$.
\end{definition}

Somewhat descriptively%
\footnote{The diagram as depicted omits all higher coherence data of the $\infty$-functor $\cC_+ p \colon N\bbDelta_+^\opp \to \scC$.},
we can depict the augmented simplicial object $\cC_+ p = \jmath_*\{p\}$ as the diagram
\begin{equation}
\begin{tikzcd}
	X &
	Y \ar[r] \ar[l]
	& Y \underset{X}{\times} Y
	\ar[l, shift left = 0.1cm] \ar[l, shift left = -0.1cm]
	\ar[r, shift left = 0.1cm] \ar[r, shift left = -0.1cm]
	& Y \underset{X}{\times} Y \underset{X}{\times} Y \ \cdots\,.
	\ar[l] \ar[l, shift left = 0.2cm] \ar[l, shift left = -0.2cm]
\end{tikzcd}
\end{equation}

\begin{definition}
We say that an $\infty$-category $\scC$ \textit{admits geometric realisations} if it has all colimits indexed by $N\bbDelta^\opp$.
Given an $\infty$-functor $D \colon N\bbDelta^\opp \to \scC$, we write
\begin{equation}
	|D| \coloneqq \colim (D \colon N\bbDelta^\opp \longrightarrow \scC)\,.
\end{equation}
\end{definition}

\begin{definition}
Let $\scC$ be an $\infty$-category with finite limits and geometric realisations.
The \textit{1-image} of a morphism $p \colon P \to X$ in $\scC$ is the colimit
\begin{equation}
	\im_1(p) \coloneqq |\cC p|\,.
\end{equation}
The \textit{1-image factorisation} of a morphism $p \colon P \to X$ is the canonical factorisation
\begin{equation}
\begin{tikzcd}
	P \ar[rr, "p"] \ar[dr]
	& & X
	\\
	& \im_1(p) \ar[ur, hookrightarrow, "\iota"']
\end{tikzcd}
\end{equation}
arising from the universal property of $|\cC p|$.
\end{definition}

The replacement for a `surjective' morphism of topological spaces is the following type of morphism:

\begin{definition}
Let $\scC$ be an $\infty$-category with finite limits and geometric realisations.
A morphism $p \colon P \to X$ is an \textit{effective epimorphism} if the canonical morphism
\begin{equation}
	\im_1(p) = |\cC p| \longrightarrow X
\end{equation}
is an equivalence in $\scC$.
\end{definition}

\begin{example}
A morphism $p \colon P \to X$ in the $\infty$-category $\scS$ is an effective epimorphism if and only if $\pi_0 p \colon \pi_0 P \to \pi_0 X$ is a bijection~\cite[Cor.~7.2.1.15]{Lurie:HTT}.
\qen
\end{example}

Effective epimorphisms can be interpreted as giving a notion of covering morphism in an $\infty$-category.
This motivates the following definition of $\infty$-bundles, generalising the notion of a locally trivial morphism with a fixed fibre:

\begin{definition}
\label{def:oo-bundle}
\cite[Def.~4.1]{NSS:oo-bundles_I}
Let $X, V$ be two objects in an $\infty$-category $\scC$ with finite limits and geometric realisations.
An \emph{$\infty$-bundle with fibre $V$ on $X$} is a morphism $p \colon E \to X$ such that there exists a pullback square
\begin{equation}
\begin{tikzcd}
	Y \times V \ar[r] \ar[d, "\pr"']
	& E \ar[d, "p"]
	\\
	Y \ar[r, "q"']
	& X
\end{tikzcd}
\end{equation}
where the bottom morphism $q \colon Y \to X$ is an effective epimorphism.
The $\infty$-category $\Bun_V(X)$ of $\infty$-bundles on $X$ with fibre $V$ is the full $\infty$-subcategory of $\scC_{/X}$ on the $\infty$-bundles on $X$ with fibre $V$.
\end{definition}

If we further pull back to the fibre product $Y \times_X Y = \cC_1 p$, we obtain a canonical equivalence
\begin{equation}
	\cC_1 p \times V \simeq (q \circ d_1)^*E \simeq (q \circ d_0)^*E \simeq \cC_1 p \times V
\end{equation}
\textit{in the slice $\infty$-category \smash{$\scC_{/\cC_1 p}$}}.
We can view this as a fibre-preserving equivalence, as we would expect for gluing a bundle with fibre $V$ on the base object $E$.


\section{From groups to higher groups}
\label{sec:Groups and actions}


Next, we need a homotopy coherent formulation of the concept of a group.
A direct way to achieve this is by using that there is a definition of groupoid objects internal to any $\infty$-category $\scC$ with pullbacks, and then viewing a group as a groupoid with only one object.
We begin with the analogous definitions for category and monoid objects internal to $\infty$-categories with pullbacks and, in a second step, implement invertibility properties to pass to groupoid and group objects, respectively.

\begin{definition}
\label{def:Cat and Mon object}
Let $\scC$ be an $\infty$-category with pullbacks.
\begin{myenumerate}
\item A \textit{category object} in $\scC$ is a simplicial object $C \colon N\bbDelta^\opp \to \scC$ satisfying the \textit{Segal conditions}:
for each $m \in \NN$, $m \geq 2$, and each partition of $[m]$ into two \textit{ordered} sets, $[r] \cup [s] = [m]$ such that the images of $[r]$ and $[s]$ in $[m]$ intersect precisely in the final object of $[r]$ and the initial object of $[s]$, the canonical morphism
\begin{equation}
	C_m \longrightarrow C_r \underset{C_0}{\times} C_s
\end{equation}
is an equivalence in $\scC$.

\item Suppose that $\scC$ also has a final object.
A \textit{monoid object} in $\scC$ is a category object $M \colon N\bbDelta^\opp \to \scC$ such that $M_0$ is a final object of $\scC$.
\end{myenumerate}
\end{definition}

Given a category object $C$ in $\scC$, we can interpret the object $C_k \in \scC$ as encoding composable $k$-tuples of 1-morphisms in $C$, together with coherent choices of compositions of these $k$-tuples.
The Segal conditions in Definition~\ref{def:Cat and Mon object} then allow us to find essentially unique compositions for any composable pair of an $r$-tuple and an $s$-tuple of morphisms, each with chosen compositions in $C$.
The key property distinguishing a groupoid and a category in this picture is that in a groupoid we can also compose $k$-tuples of morphisms where some (or all) of the morphisms face the wrong direction.

To make this formal, we first need the following convention:
given a simplicial object $X \colon N\bbDelta^\opp \to \scC$ in an $\infty$-category $\scC$ and a simplicial set $K \in \sSet$, we set (if the following limits exists in $\scC$)
\begin{equation}
	X(K) \coloneqq \lim \big( N(\bbDelta_{/K})^\opp \to N\bbDelta^\opp \to \scC \big)\,,
\end{equation}
where the first arrow arises from the projection $\bbDelta_{/K} \to \bbDelta$ and the second is the functor $X$.

\begin{remark}
The above limit is guaranteed to exist whenever $K$ is finite and $\scC$ is finitely complete, or for each small simplicial set $K \in \sSet$ if $\scC$ is complete.
It is the value of the right Kan extension of $X$ along the Yoneda embedding $\bbDelta \hookrightarrow \sSet$.
\qen
\end{remark}

Any subset $R \subset [m]$, for $m \in \NN$, defines a canonical inclusion $\Delta^{|R|} \hookrightarrow \Delta^m$.

\begin{definition}
\label{def:Gpd and Grp object}
\cite[Def.~6.1.2.7, Prop.~6.1.2.6]{Lurie:HTT}
Let $\scC$ be an $\infty$-category with pullbacks.
\begin{myenumerate}
\item A \textit{groupoid object} in $\scC$ is a simplicial object $C \colon N\bbDelta^\opp \to \scC$ satisfying the \textit{groupoidal Segal conditions}:
for each $m \in \NN$, $m \geq 2$, and each partition of $[m]$ into two (unordered) sets, $R \cup S = [m]$ such that the images of $R$ and $S$ in $[m]$ intersect in precisely one object of $[m]$, the canonical morphism
\begin{equation}
	C_m \longrightarrow C(\Delta^{|R|}) \underset{C_0}{\times} C(\Delta^{|S|})
\end{equation}
is an equivalence in $\scC$.

\item Suppose that $\scC$ also has a final object.
A \textit{group object} in $\scC$ is a groupoid object $G \colon N\bbDelta^\opp \to \scC$ such that $G_0$ is a final object of $\scC$.
\end{myenumerate}
\end{definition}

Equivalently, a groupoid object in $\scC$ is a category object $C \in \scC$ such that the morphism induced by the inclusion $\Lambda^2_0 \hookrightarrow \Delta^2$ is an equivalence
\begin{equation}
	X(\Delta^2) \longrightarrow X(\Lambda^2_0)\,.
\end{equation}
This can be found in the proof of~\cite[Prop.~1.1.8]{Lurie:Goodwillie}.

We write $\Gpd(\scC)$ and $\Grp(\scC)$ for the full $\infty$-subcategories of $\scFun(N\bbDelta^\opp, \scC)$ on the groupoid and group objects in $\scC$, respectively.
With the notion of groupoid objects at hand, we can now define the type of $\infty$-category which is the natural home for principal $\infty$-bundles:

\begin{definition}
\label{def:oo-topos}
\cite[Def.~6.1.0.4, Thm.~6.1.0.6]{Lurie:HTT}
An \textit{$\infty$-topos} is an $\infty$-category $\scX$ satisfying the following properties:
\begin{myenumerate}
\item \textit{(Presentability)}
The $\infty$-category $\scX$ is presentable%
\footnote{See~\cite[Def.~5.5.0.1]{Lurie:HTT} for a definition of presentability.
Equivalently, an $\infty$-category is presentable if it is a reflective localisation of a presheaf $\infty$-category at a small class of morphisms~\cite[Def.~7.11.15]{Cisinski:HCats_HoAlg}.}.
In particular, it has all small limits and colimits~\cite[Cor.~5.5.2.4]{Lurie:HTT}.
We denote its initial object by $\emptyset \in \scX$ and its final object by $* \in \scX$.

\item \textit{(Colimits are universal in $\scX$)}
For each morphism $f \colon A \to B$ in $\scX$ and each diagram $D \colon I \to \scX_{/B}$, the canonical morphism
\begin{equation}
	\underset{i \in I}{\colim} \big( A \underset{B}{\times} Di \big)
	\longrightarrow A \underset{B}{\times} \underset{i \in I}{\colim}\, Di
\end{equation}
is an equivalence in $\scX$.

\item \textit{(Coproducts are disjoint in $\scX$)}
For each $A, B \in \scX$, the pushout square
\begin{equation}
\begin{tikzcd}
	\emptyset \ar[r] \ar[d]
	& A \ar[d]
	\\
	B \ar[r]
	& A \sqcup B
\end{tikzcd}
\end{equation}
is also a pullback square.

\item \textit{(Groupoids are effective in $\scX$)}
That is, given a groupoid object $C \colon N\bbDelta^\opp \to \scX$, consider the extension $C_+ \colon (N\bbDelta^\opp)^\triangleright \cong N\bbDelta_+^\opp \longrightarrow \scX$ obtained by forming the colimit of $C$.
Using the notation of Definition~\ref{def:Cech nerves}, the identity $\jmath^*C_+ = (p \colon C_0 \to |C|)$ together with the adjunction $\jmath^* \dashv \jmath_*$ induce a canonical morphism $C_+ \longrightarrow \jmath_*\{p\} = \cC_+ p$.
The condition is that this morphism is an equivalence of augmented simplicial objects $N\bbDelta_+^\opp \to \scX$.
In particular, it then follows that $p \colon C_0 \to |C|$ is an effective epimorphism.
\end{myenumerate}
\end{definition}

\begin{example}
Important examples of $\infty$-topoi consist of the $\infty$-category $\scS$ of spaces, $\infty$-presheaves $\PSh(\scC)$ of spaces on any $\infty$-category $\scC$, and any accessible, left exact, reflective localisation of $\infty$-categories of the form $\PSh(\scC)$~\cite[Prop.~6.1.5.3]{Lurie:HTT}.
In particular, each sheaf $\infty$-category on an $\infty$-site $(\scC, \tau)$ (or, equivalently, each topological localisation of a presheaf $\infty$-category $\PSh(\scC)$) is an $\infty$-topos~\cite[Cor.~6.2.1.7, Prop.~6.2.2.9]{Lurie:HTT}.
\qen
\end{example}

In each $\infty$-topos $\scX$ there are reflective localisations~\cite[p.~587]{Lurie:HTT}
\begin{equation}
\begin{tikzcd}
	\scFun(N\bbDelta^\opp, \scX) \ar[r, shift left=0.15cm, "\perp"' yshift=0.05cm]
	& \Gpd(\scX) \ar[r, shift left=0.15cm, "\perp"' yshift=0.05cm] \ar[l, hookrightarrow, shift left=0.15cm]
	& \Grp(\scX)\,. \ar[l, hookrightarrow, shift left=0.15cm]
\end{tikzcd}
\end{equation}
Furthermore, there is a canonical equivalence~\cite[Lemma~7.2.2.11]{Lurie:HTT}
\begin{equation}
\label{eq:Omega--B equivalence}
\begin{tikzcd}
	\Omega : \scX^{*/}_{\geq 1} \ar[r, shift left=0.15cm, "\perp"' yshift=0.05cm]
	& \Grp(\scX) : \rmB \ar[l,shift left=0.15cm]
\end{tikzcd}
\end{equation}
between the $\infty$-category of pointed, connected objects in $\scX$ and the $\infty$-category of group objects in $\scX$.
The functor $\Omega$ sends a pointed connected object $(* \to X)$ to the \v{C}ech nerve $\cC (* \to X)$, and the functor $\rmB$ takes a group object $G$ to the colimit $|G|$ of its underlying groupoid object.

\begin{remark}
\label{rmk:presenting oo-groups}
For the $\infty$-topos $\scX = \scS$, the equivalence~\eqref{eq:Omega--B equivalence} has a presentation in terms of localisations of 1-categories via McLane's simplicial delooping functor:
this establishes a Quillen equivalence between model categories of simplicial groups and reduced simplicial sets~\cite[Cor.~3.34]{NSS_oo-Bundles_II}.
The presentation further carries over to sheaf $\infty$-topoi which satisfy a mild condition~\cite[Prop.~3.35]{NSS_oo-Bundles_II}.
\qen
\end{remark}

\begin{remark}
Monoid objects and group objects in $\infty$-topoi can be defined in terms of algebras over the associative $\infty$-operad $\bbE_1$; we comment on this in more detail in Remark~\ref{rmk:Grp objs via E_n-algebras} below.
\qen
\end{remark}


\section{Principal $\infty$-bundles}
\label{sec:Pr-oo-bundles}


In this section we take the final two steps in building up definition of principal $\infty$-bundles in $\infty$-topoi.
We define group actions in $\infty$-topoi and and formulate the $\infty$-categorical version of the principality condition.
For the following definition, see~\cite[Def.~3.1]{NSS:oo-bundles_I}, \cite[Def.~3.15]{Bunk:Pr-oo-Bundles}, and~\cite[Def.~13.1.25]{ADH:Diff_Coho} for a slightly different, but equivalent, formulation.
See also~\cite[Chs.~4, 5]{Lurie:HA} and~\cite[Sec.~13]{ADH:Diff_Coho} for more on group objects and group actions in $\infty$-categories.

\begin{definition}
\label{def:group-action}
Let $\scC$ be an $\infty$-category with finite products and $G$ a group object in $\scC$.
A \textit{$G$-action on an object $P \in \scC$} (from the right) is a simplicial object $P \dslash G$ in $\scC$ satisfying that
\begin{myenumerate}
\item $(P \dslash G)_n = P \times G_1^n$, for each $n \in \NN_0$

\item the face map $d_1 \colon P \times G_1 \to P$ coincides with the canonical projection onto $P$,

\item the degeneracy map $s_0 \colon P \to P \times G_1$ coincides with the product of the identity $\id_P$ with the degeneracy map $s_0 \colon * \to G_1$ of the simplicial object $G$, and

\item the projection morphisms $P \times G_1^n \to G_1^n$ define a morphism $P \dslash G \to G$ of simplicial objects in $\scC$.
\end{myenumerate}
A \textit{morphism of $G$-actions} is a morphism $f \colon P \dslash G \to Q \dslash G$ of simplicial objects in $\scC$ such that the following triangle in $\scFun(N\bbDelta^\opp, \scX)$ commutes:
\begin{equation}
\begin{tikzcd}[column sep=0.5cm]
	P \dslash G \ar[rr, "f"] \ar[dr]
	& & Q \dslash G \ar[dl]
	\\
	& G &
\end{tikzcd}
\end{equation}
\end{definition}

\begin{remark}
In~\cite[Def.~4.2.2.2]{Lurie:HA}, Lurie gives a definition of a left action of a monoid object in an $\infty$-category $\scC$, and upon replacing $\{n\}$ by $\{0\}$ in axiom (ii) of that definition, we obtain a definition of a right action.
On first inspection, Definition~\ref{def:group-action} appears stricter than~\cite[Def.~4.2.2.2]{Lurie:HA}; however, each action in the latter sense is equivalent to one in the sense of Definition~\ref{def:group-action} above by~\cite[Lemma~3.12]{Bunk:Pr-oo-Bundles}.
\qen
\end{remark}

We can think of the simplicial object $P \dslash G$ as the action groupoid associated to a $G$-action on $P$:

\begin{theorem}
\label{st:group actions are groupoid objects}
\emph{\cite[Thm.~3.19]{Bunk:Pr-oo-Bundles}}
Let $\scC$ be an $\infty$-category with finite limits.
For each group object $G \in \Grp(\scC)$ and each $G$-action $P \dslash G$ in $\scC$, the simplicial object $P \dslash G$ is a groupoid object in $\scC$.
\end{theorem}

\begin{definition}
In the above setting, let $p \colon P \to X$ be a morphism in $\scC$.
A \textit{$G$-action on $P$ over $X$} is an augmented simplicial object $(P \dslash G \to X) \in \scFun(N\bbDelta_+^\opp, \scC)$ such that
\begin{myenumerate}
\item the restriction of the functor $(P \dslash G \to X) \in \scFun(N\bbDelta_+^\opp, \scC)$ to $N\bbDelta^\opp$ is a $G$-action $P \dslash G$ on $P$, and

\item we have that $(P \dslash G \to X)_{-1} = X$.
\end{myenumerate}
A \textit{morphism of $G$-actions over $X$} is a morphism $(P \dslash G \to X) \longrightarrow (Q \dslash G \to X)$ in $\scFun(N\bbDelta_+^\opp, \scC)$ whose underlying morphism of simplicial objects is a morphism of $G$-actions, and whose component in degree $-1$ is the identity.
\end{definition}

\begin{example}
\label{eg:oo-group actions}
Let $\scX$ be an $\infty$-topos.
There are various important examples of group actions:
\begin{myenumerate}
\item For any group object $G$ in $\scX$, the simplicial object $G \in \scFun(N\bbDelta^\opp, \scX)$ canonically encodes an action $* \dslash G$ of $G$ on the final object $* \in \scX$.
We will therefore use the notation $G$ and $* \dslash G$ interchangeably.
Given any group action $P \dslash G$ in $\scX$, the collapse morphisms $P \to *$ induces a canonical morphism of $G$-actions $P \dslash G \to * \dslash G$.
Further, note that $\rmB G = |{*} \dslash G|$.

\item The decalage of the simplicial object $G = * \dslash G$ yields a $G$ action $G_1 \dslash G$; this is the canonical action of the group object $G \in \Grp(\scX)$ on its underlying object $G_1 \in \scX$ via right multiplication (see~\cite[Ex.~3.13]{Bunk:Pr-oo-Bundles} for details).

\item Let $x_0 \colon * \to \scX$ be a pointed object in $\scX$.
We obtain a \textit{based loop group} $\Omega_{x_0} X \coloneqq \cC x_0 \in \Grp(\scX)$.
For $\scX = \scS$, this is the well-known grouplike $\bbE_1$-structure on $\Omega_x X$ (compare Section~\ref{sec:Motivation}).
If we define $P_{x_0}X \in \scX$ as the pullback
\begin{equation}
\begin{tikzcd}
	P_{x_0}X \ar[r] \ar[d]
	& X \ar[d, equal]
	\\
	* \ar[r, "\{x_0\}"']
	& X
\end{tikzcd}
\end{equation}
we can interpret this as the based path-space object of $(X, x_0)$.
Indeed, if $\scX = \scS$, there is a canonical equivalence \smash{$P_{x_0} X \simeq X_{x_0/}$} of $\infty$-groupoids.
The \v{C}ech nerve of the canonical equivalence $p \colon P_{x_0} X \to X$ comes with a canonical equivalence $(\cC p)_n \simeq P_{x_0}X \times (\Omega_{x_0} X)^n$, for each $n \in \NN_0$, which exhibits a group action $P_{x_0} X \dslash \Omega_{x_0} X \longrightarrow X$ in $\scX$ over $X$ (using~\cite[Lemma~3.12]{Bunk:Pr-oo-Bundles}).

\item Given any group action $P \dslash G$ in $\scX$, we obtain an augmented simplicial object $(P \dslash G \to |P \dslash G|)$ in $\scX$ by appending the colimit of the simplicial object $P \dslash G$.
This is a $G$-action on $P$ over $|P \dslash G|$.

\item Each morphism $A \to G$ in $\Grp(\scX)$ induces a canonical action of $A$ on the underlying object $G_1 \in \scX$.
We can think of this as mapping elements of $A$ to $G$ and then acting by means of the multiplication in $G$~\cite[Prop.~3.24]{Bunk:Pr-oo-Bundles}.
\qen
\end{myenumerate}
\end{example}

The following definition was introduced in~\cite{NSS:oo-bundles_I}:

\begin{definition}
\label{def:G-pr oo-Bun}
\cite[Def.~3.4]{NSS:oo-bundles_I}
Let $\scX$ be an $\infty$-topos and $G \in\Grp(\scX)$ a group object in $\scX$.
A \textit{$G$-principal $\infty$-bundle} over $X \in \scX$ is a $G$-action $(P \dslash G \to X)$ over $X$ whose underlying augmented simplicial object $(P \dslash G \to X) \in \scFun(N\bbDelta_+^\opp, \scX)$ is a colimiting cocone for the simplicial diagram $P \dslash G \in \scFun(N\bbDelta^\opp, \scX)$.
Equivalently, the canonical morphism $|P \dslash G| \to X$ is an equivalence in $\scX$.

A \textit{morphism of $G$-principal $\infty$-bundles} $(P \dslash G \to X) \longrightarrow (Q \dslash G \to X)$ is a morphism of the underlying $G$-actions over $X$.
This defines, for each $X \in \scX$, a space $\Bun(X;G) \in \scS$.
\end{definition}

\begin{remark}
A priori, $\Bun(X;G)$ is an $\infty$-category, rather than an $\infty$-groupoid.
It is a non-trivial result that, in fact, every morphism of $G$-principal $\infty$-bundles over any object $X \in \scX$ is an equivalence (see Theorem~\ref{st:Class Thm for G-pr-oo-Bdls} and Corollary~\ref{st:Bun(-;G) is space-valued} below).
This is the $\infty$-categorical analogue of Remark~\ref{rmk:classical pr-Bun mps are isos}.
\qen
\end{remark}

\begin{definition}
\label{def:shear mp}
Let $G$ be a group object in an $\infty$-topos $\scX$, and let $(P \dslash G \to X)$ be a $G$-action over $X$.
We set $p \coloneqq d_1 \colon P \dslash G \to X$.
By the procedure in Definition~\ref{def:oo-topos}(4) we obtain a canonical morphism $P \dslash G \to \cC p$ in $\scFun(N\bbDelta^\opp, \scX)$.
We call this the \textit{shear morphism} of the $G$-action $P \dslash G \to X$.
\end{definition}

\begin{remark}
It follows directly from Definition~\ref{def:G-pr oo-Bun} that, for each $G$-principal $\infty$-bundle $(P \dslash G \to X)$, the augmenting morphism $p \coloneqq d_{-1} \colon P \to X$ is an effective epimorphism in $\scX$:
by Theorem~\ref{st:group actions are groupoid objects} we know that $P \dslash G$ is a groupoid object in $\scX$.
Since groupoid objects in $\scX$ are effective (Definition~\ref{def:oo-topos}(4)) the shear morphism $P \dslash G \to \cC p$ is an equivalence in $\Gpd(\scX)$.
It then follows that $p$ is indeed an effective epimorphism.
\qen
\end{remark}

We can characterise $G$-principal $\infty$-bundles in an alternative way, which more closely resembles the classical Definition~\ref{def:ord principal bundle}.

\begin{definition}
\label{def:principality}
Let $G$ be a group object in an $\infty$-topos $\scX$.
A $G$-action $(P \dslash G \to X)$ over $X$ is called \textit{principal} (or a \textit{torsor}) if the following equivalent conditions are satisfied:
\begin{myenumerate}
\item The shear morphism $P \dslash G \to \cC p$ is an equivalence in $\Gpd(\scX)$.

\item The canonical diagram
\begin{equation}
\begin{tikzcd}[column sep=1.25cm, row sep=1cm]
	P \times G_1 \ar[r, "d_1 = \pr_P"] \ar[d, "d_0"']
	& P \ar[d, "p"]
	\\
	P \ar[r, "p"']
	& X
\end{tikzcd}
\end{equation}
is a pullback diagram.
\end{myenumerate}
\end{definition}

The equivalence of the conditions in Definition~\ref{def:principality} is shown in~\cite[Prop.~3.29]{Bunk:Pr-oo-Bundles}.

\begin{proposition}
\label{st:alt char of pr-oo-Buns}
\emph{\cite[Prop.~3.31]{Bunk:Pr-oo-Bundles}}
Let $\scX$ be an $\infty$-topos and $G \in\Grp(\scX)$ be a group object in $\scX$.
A $G$-principal $\infty$-bundle over $X \in \scX$ is equivalently a $G$-action $P \dslash G \longrightarrow X$ over $X$ such that
\begin{myenumerate}
\item the underlying morphism $P \to X$ is an effective epimorphism and

\item the $G$-action is principal (in the sense of Definition~\ref{def:principality}).
\end{myenumerate}

\end{proposition}

\begin{example}
\label{eg:pr-oo-bundles}
Let $\scX$ be an $\infty$-topos and $G \in \Grp(\scX)$ a  group object.
\begin{myenumerate}
\item The augmented simplicial object $* \dslash G \to BG = |{*} \dslash G|$ is a $G$-principal $\infty$-bundle in $\scX$.
It is called the \textit{universal $G$-principal $\infty$-bundle}.
It generalises the classical universal bundle $\rmE G \to \rmB G$ from algebraic topology, for a topological group $G$.
Note that in this classical setting $\rmE G$ is a contractible space with a free $G$-action.
In the present, $\infty$-categorical setting we are free to replace these data by the trivial $G$-action on the final object $* \in \scX$.

\item The augmented simplicial object $G_1 \dslash G \to |G_1 \dslash G| \simeq *$ from Example~\ref{eg:oo-group actions}(2) is a $G$-principal $\infty$-bundle.
It is the trivial principal $G$-bundle over the final object $* \in \scX$.

\item Given any object $X \in \scX$, then $(X \times G_1 \dslash G) \to X$ is the \textit{trivial} $G$-principal $\infty$-bundle on $X$.
A $G$-principal $\infty$-bundle over $X$ which is equivalent to the trivial bundle is called \textit{trivialisable}.

\item Let $x_0 \colon * \to X$ be a pointed object in $\scX$.
The $\Omega_{x_0} X$-action $P_{x_0} X \dslash \Omega_{x_0} X \to X$ over $X$ in Example~\ref{eg:oo-group actions}(3) is an $\Omega_{x_0} X$-principal $\infty$-bundle if and only if $X$ is connected.
Here the critical condition is that the morphism $P_{x_0} X \to X$ needs to be an effective epimorphism.
In particular, for $\scX = \scS$ this establishes the introductory example from Section~\ref{sec:Motivation} as a principal $\infty$-bundle.
\qen
\end{myenumerate}
\end{example}

\begin{remark}
Ordinary principal bundles, in the sense of Definition~\ref{def:ord principal bundle}, are, in particular, $G$-principal $\infty$-bundles in the $\infty$-category $N\Top$.
\qen
\end{remark}

Principal $\infty$-bundles have applications in the higher algebra of groups in $\infty$-topoi:

\begin{definition}
\label{def:oo-group extension}
\cite[Def.~4.26]{NSS:oo-bundles_I}
A pair of composable morphisms $A \to G \to H$ in $\Grp(\scX)$ is an \emph{extension of $\infty$-groups} if the induced sequence $\rmB A \to \rmB G \to \rmB H$ is a fibre sequence in $\scX$.
\end{definition}

\begin{theorem}
\emph{\cite[Thm.~3.48]{Bunk:Pr-oo-Bundles}}
A pair of composable morphisms $A \to G \to H$ in $\Grp(\scX)$ is an extension of $\infty$-groups if and only if the morphism $G_1 \to H_1$ on underlying objects in $\scX$, together with the induced action of $A$ on $H_1$ (see Example~\ref{eg:oo-group actions}(5)) is an $A$-principal $\infty$-bundle in $\scX$.
\end{theorem}


\section{Non-abelian cohomology and classifying objects}
\label{sec:Non-ab coho and BG}


Let $\scX$ be an $\infty$-topos and $G \in \Grp(\scX)$.
One can show that $G$-principal $\infty$-bundles can be pulled back along morphisms  in $\scX$:

\begin{proposition}
\emph{\cite[Props.~3.33, 3.41\footnote{The functoriality is not explicitly stated there, but the proof is fully functorial.}]{Bunk:Pr-oo-Bundles}}
Let $f \in \scX(X,Y)$ and $(P \dslash G \to Y) \in \Bun(Y;G)$.
Then, the pullback $f^*P \coloneqq X \times_Y P$ carries a natural $G$-action over $X$, and this makes $(f^*P) \dslash G \to X$ into a $G$-principal $\infty$-bundle, denoted $f^*(P \dslash G \to Y)$.
We obtain a morphism of spaces (see Theorem~\ref{st:Class Thm for G-pr-oo-Bdls} and Corollary~\ref{st:Bun(-;G) is space-valued} below)
\begin{equation}
	\scX(X,Y) \longrightarrow \scS \big( \Bun(Y;G), \Bun(X;G) \big)\,.
\end{equation}
Further, each $G$-principal $\infty$-bundle arises in an essentially unique way as a pullback of the universal $G$-principal $\infty$-bundle $* \dslash G \to \rmB G$.
\end{proposition}

The following classification theorem for principal $\infty$-bundles in an $\infty$-topos $\scX$ is~\cite[Prop.~3.13, Thm.~3.17]{NSS:oo-bundles_I}, (see also \cite[Props.~3.33, 3.41]{Bunk:Pr-oo-Bundles} for a more detailed treatment of the essential-surjectivity part of the statement).

\begin{theorem}
\label{st:Class Thm for G-pr-oo-Bdls}
For each $X \in \scX$, the pullback of $G$-principal $\infty$-bundles induces an equivalence
\begin{equation}
	\scX(X, \rmB G) \longrightarrow \Bun(X;G)\,,
	\qquad
	f \longmapsto f^*(* \dslash G \to \rmB G)\,.
\end{equation}
\end{theorem}

\begin{corollary}
\label{st:Bun(-;G) is space-valued}
The assignment $(X \in \scX) \longmapsto \Bun(X;G)$ is a functor $\Bun(-;G) \colon \scX^\opp \to \scS$, which classifies the canonical right fibration \smash{$\scX_{/\rmB G} \to \scX$} in $\sSet$.
In particular, each morphism of $G$-principal $\infty$-bundles is an equivalence.
\end{corollary}

\begin{remark}
Recall from Remark~\ref{rmk:Grp objs via E_n-algebras} that $\infty$-groups in sufficiently nice $\infty$-topoi have presentations in terms of simplicial groups.
Similarly, principal $\infty$-bundles in such $\infty$-topoi can be presented by means of 1-categorical constructions in simplicial homotopy theory (see, in particular, the results involving weakly principal bundles for simplicial groups~\cite[Def.~3.79, Thm.~3.95]{NSS_oo-Bundles_II}).
\qen
\end{remark}

\begin{remark}
Given \textit{any} $G$-action $V \dslash G$ in $\scX$, the canonical morphism $V \dslash G \to |V \dslash G|$ is a $G$-principal $\infty$-bundle.
Indeed, $V \dslash G$ is a groupoid object in $\scX$, so that the morphism $p \colon V \to |V \dslash G|$ is an effective epimorphism.
In particular, we obtain that the $\infty$-category of $G$-actions is also canonically equivalent to the overcategory $\scX_{/\rmB G}$.
\qen
\end{remark}

\begin{example}
\label{eg:pr-oo-bdls from classifying objects}
The description of $G$-principal $\infty$-bundles via the classifying object $\rmB G$ is often useful in practise:
\begin{myenumerate}
\item Recall the $\infty$-categories $\scH$ and $\scH_\tau$ from Example~\ref{eg:oo-topoi of smooth spaces}.
There is a fully faithful embedding $\Mfd \hookrightarrow \scH_\tau$ of the category of smooth manifolds and smooth maps into $\scH_\tau$.
This sends the abelian Lie group $\rmU(1)$ to an $\bbE_\infty$-group object in $\scH_\tau$.
Thus, there exists a delooping $\rmB^n \rmU(1)$, for each $n \in \NN$.
Under the presentation in Example~\ref{eg:oo-topoi of smooth spaces} this corresponds to the simplicial presheaf obtained by applying the Dold-Kan correspondence to the homotopy sheaf $\rmU(1)[n]$ of chain complexes of abelian groups.
The \textit{$\infty$-groupoid of $n$-gerbes}, also called \textit{$\rmU(1)$-$(n{+}1)$-bundles}, or \textit{$\rmB^n\rmU(1)$-principal $\infty$-bundles}, on an object $X \in \scH_\tau$ is the mapping space
\begin{equation}
	\Grb^n(X) \coloneqq \Bun \big( X; \rmB^n\rmU(1) \big)
	= \scH_\tau \big( X, \rmB^{n+1} \rmU(1) \big)\,.
\end{equation}
Under the presentation of $\scH_\tau$ from Example~\ref{eg:oo-topoi of smooth spaces}, this space can be modelled by the simplicial hom space in $\Fun(\Cart^\opp, \sSet)$ from a cofibrant object presenting $X$ in the \textit{$\tau$-local projective} model structure.
For instance, if $X$ is the image of a manifold $M$ under the embedding $\Mfd \hookrightarrow \scH_\tau$, and $\CU = \{U_a\}_{a \in \Lambda}$ is a good open cover of $M$, then the \v{C}ech nerve of the cover is a cofibrant object as desired (see, for instance,~\cite{FSS:Cech_diff_char_classes_via_L_infty, Schreiber:DCCT, Bunk:Higher_Sheaves} for more on this).

\item For each $n \in \NN$, there also exists an object $\rmB^n_\nabla \rmU(1) \in \scH_\tau$ which classifies $(n{-}1)$-gerbes on $M$ \textit{with connection}.
It is presented by the simplicial homotopy sheaf obtained via the Dold-Kan correspondence from the Deligne complex of sheaves of abelian groups~\cite{FSS:Cech_diff_char_classes_via_L_infty, Schreiber:DCCT}
\begin{equation}
\begin{tikzcd}
	\rmU(1) \ar[r, "\dd \log"]
	& \Omega^1 \ar[r, "\dd"]
	& \cdots \ar[r, "\dd"]
	& \Omega^n\,.
\end{tikzcd}
\end{equation}
Each of these objects is again an abelian group object, and so admits deloopings $\rmB^k \rmB^n_\nabla \rmU(1)$ for each $k \in \NN$.
For instance, $\rmB_\nabla \rmU(1)$-principal $\infty$-bundles (equivalently known as 1-gerbes with connective structure) on manifolds are closely related to exact Courant algebroids via Hitchin's generalised tangent bundle construction~\cite[Sec.~16]{BS:Higher_Syms_and_Deligne_Coho}.

\item Given a ring spectrum $R$ in the $\infty$-category of spaces, one obtains a group object $\mathrm{GL}_1(R) \in \Grp(\scS)$ of its units.
Given a space $X \in \scS$, one can interpret maps $X \to \rmB \mathrm{GL}_1(R)$ as $R$-line bundles on $X$ (with a flat connection); these objects govern the twisted $R$-(co)homology of $X$~\cite{ABGHR} (see also~\cite{DY:Adams_spec_seq_and_Thom_sp} for an overview).
\qen
\end{myenumerate}
\end{example}

Theorem~\ref{st:Class Thm for G-pr-oo-Bdls} also implies that $G$-principal $\infty$-bundles are cocycles for non-abelian cohomology in $\scX$:

\begin{definition}
\cite[Def.~2.24]{NSS:oo-bundles_I}
Given an object $T \in \scX$, we define, for each $X \in \scX$, the \textit{cohomology set of $X$ with coefficients in $T$} as
\begin{equation}
	\rmH^0(X;T) \coloneqq \pi_0 \scX(X,T)\,.
\end{equation}
For $n \in \NN$, we can always define the \textit{$(-n)$-th cohomology group of $X$ with coefficients in $T$} as
\begin{equation}
	\rmH^{-n}(X;T) \coloneqq \rmH^0(X; \Omega^n T)\,,
\end{equation}
If $T \in \scX_{\geq 1}^{*/}$ is connected and pointed, the equivalence~\eqref{eq:Omega--B equivalence} provides a canonical equivalence
\begin{equation}
	\rmH^0(X;T) = \pi_0 \scX(X,T)
	\simeq \pi_0 \scX(X, \rmB \Omega T)
	\simeq \pi_0 \Bun(X; \Omega T)\,.
\end{equation}
\end{definition}

That is, $\Omega T$-principal $\infty$-bundles are cocycles for cohomology with coefficients in $T$.
Equivalently, cohomology with coefficients in $T$ classifies $\Omega T$-principal $\infty$-bundles.

\begin{definition}
\cite[Def.~2.24 (ctd.)]{NSS:oo-bundles_I}
If $T$ is an $n$-fold loop object, i.e.~there exists a \textit{$k$-fold delooping} $\rmB^k T \in \scX^{*/}$, for $k = 1, \ldots, n$, such that $\rmB^{k-1} T \simeq \Omega \rmB^k T$ as group objects in $\scX$, then we can also define the \textit{$n$-th cohomology set of $X$ with coefficients in $T$} as
\begin{equation}
	\rmH^n(X;T) \coloneqq \rmH^0(X; \rmB^n T)
	= \pi_0 \scX(X, \rmB^n T)
	\simeq \pi_0 \Bun(X; \rmB^{n-1}T)\,.
\end{equation}
\end{definition}

In this case, for $1 \leq k \leq n$, we also find that
\begin{equation}
	\rmH^{n-k}(X;T)
	= \pi_0 \scX(X, \rmB^{n-k} T)
	\simeq \pi_0 \scX(X, \Omega^k \rmB^n T)
	\simeq \pi_k \scX(X, \rmB^n T)
	\simeq \pi_k \Bun(X; \rmB^{n-1}T)\,.
\end{equation}
In particular, if $T$ is an $n$-fold loop object, then $\rmH^n(X;T)$ is a set, $\rmH^{n-1}(X;T)$ is a group, and $\rmH^{n-k}(X;T)$ is an abelian group for all $k \geq 2$.
This also implies descriptions of the higher cohomology groups with coefficients in $T$ in terms of principal $\infty$-bundles arising from deloopings of $T$ and their automorphisms.

\begin{remark}
\label{rmk:Grp objs via E_n-algebras}
There is an ample supply of group objects and group objects with higher deloopings from grouplike $\bbE_k$-monoids in $\scX$:
let $\bbE_k$ denote the $\infty$-operad of $k$-dimensional cubes.
Group objects can be obtained from $\bbE_1$ algebras:
by~\cite[Rmk.~5.2.6.5, Ex.~5.2.6.13]{Lurie:HA}, $\bbE_1$-algebras in $\scX$ give rise to monoid objects, and an $\bbE_1$-algebra in $\scX$ is grouplike (see~\cite[Def.~5.2.6.2]{Lurie:HA}) if and only if its associated monoid object in $\scX$ is a group object.
Using the Dunn Additivity Theorem~\cite[5.1.2.2]{Lurie:HA} one can then enhance the equivalence~\eqref{eq:Omega--B equivalence} to an equivalence, for each $k \in \NN$,
\begin{equation}
	\scX^{*/}_{\geq k} \longrightarrow \mathrm{Mon}_{\bbE_k}^{\mathrm{gp}}(\scX)
\end{equation}
between pointed, $k$-connective objects in $\scX$ and grouplike $\bbE_k$-monoids in $\scX$~\cite[5.2.6.15]{Lurie:HA}.
\qen
\end{remark}


\section{Associated $\infty$-bundles and automorphism groups}
\label{sec:associated bundles}


A fundamental construction in the theory of classical fibre bundles is the \textit{Borel construction}, or \textit{associated bundle construction}:
given a topological group $G$, a principal $G$-bundle $P \to X$ (in the classical sense, Definition~\ref{def:ord principal bundle}) on a topological space $X$, and a left action $G \circlearrowright V$ of $G$ on some topological space $V$, we can form the quotient of $P \times V$ by the induced \textit{diagonal} $G$-action:
\begin{equation}
\label{eq:classical Borel quotient}
	P \times_G V \coloneqq (P \times V)/{\sim}\,,
	\qquad
	(x,v) \simeq (x g, g^{-1} v)\,,
	\quad
	\forall\, x \in X,\, v \in V,\, g \in G\,.
\end{equation}
This canonically exhibits the structure of a fibre bundle over $X$ with typical fibre $V$.

In order to obtain a version of this construction for $\infty$-bundles, we have to categorify it and make it internal to an arbitrary $\infty$-topos $\scX$.
To that end, let $G \in \Grp(\scX)$ be a group object, $P \dslash G \to X$ a $G$-principal $\infty$-bundle, and $V \dslash G$ a $G$-action%
\footnote{Note that in the classical construction we used a \textit{left} action of $G$ on $V$---this is how this construction is usually encountered---but in the actual Borel quotient~\eqref{eq:classical Borel quotient}, this is transformed into its associated right action by acting on $V$ with $g^{-1}$ instead of $g$.
In the $\infty$-categorical case we use this right action straight away.}
on an object $V \in \scX$.
The pullback
\begin{equation}
\begin{tikzcd}
	(P \times V) \dslash G \ar[r] \ar[d]
	& V \dslash G \ar[d]
	\\
	P \dslash G \ar[r]
	& * \dslash G
\end{tikzcd}
\end{equation}
in $\scFun(N\bbDelta^\opp, \scX)$ encodes the diagonal action of $G$ on $P \times V$~\cite[Rmk.~4.3]{NSS:oo-bundles_I}.

\begin{definition}
Let $\scX$ be an $\infty$-topos, $G \in \Grp(\scX)$, $P \dslash G \to X$ a $G$-principal $\infty$-bundle, and $V \dslash G$ a $G$-action.
We define the \textit{Borel construction}, or \textit{associated $\infty$-bundle} of the above data as the colimit
\begin{equation}
	P \times_G V \coloneqq  \big| (P \times G) \dslash G \big|
	\quad \in \scX\,.
\end{equation}
\end{definition}

This is indeed an $\infty$-bundle in the sense of Definition~\ref{def:oo-bundle} by~\cite[Prop.~4.8]{NSS:oo-bundles_I}.

We will now show that in an $\infty$-topos each $\infty$-bundle $q \colon E \to X$ arises as an associated bundle, as long  as its fibre $V$ satisfies a certain size condition.
This was already discovered in the original paper~\cite[Sec.~4.1]{NSS:oo-bundles_I}.
The key step is the realisation that the structural properties of an $\infty$-topos allow us to obtain---at least at the abstract level---the automorphism $\infty$-group of any object $V \in \scX$ as an object in $\Grp(\scX)$, together with its natural action on $V$, as we now describe.

The property of $\infty$-topoi which facilitates this is the existence of \textit{classifying objects for relatively $\kappa$-small morphisms}.
In other words, any morphism $Y \to X$ in an $\infty$-topos $\scX$ which satisfies a certain size condition (see below for details) can be written, in a unique way, as the pullback of a particular morphism, called $\Obj^\kappa_* \to \Obj^\kappa$.

\begin{remark}
One can view this as an analogue, internally to $\scX$, of the existence of a universal left fibration $\scS_* \to \scS$ which classifies $\infty$-functors valued in the $\infty$-category $\scS$ of spaces.
\qen
\end{remark}

Let us make this more precise:
let $\kappa$ be a regular cardinal.
An object $X$ in an $\infty$-category $\scC$ is called \textit{$\kappa$-compact} if the functor
\begin{equation}
	\scC(X, -) \colon \scC \to \scS
\end{equation}
corepresented by $X$ preserves $\kappa$-filtered colimits%
\footnote{Recall that an $\infty$-category $\scI$ is $\kappa$-filtered if, for each $\kappa$-small simplicial set $K$ and each morphism $K \to \scI$, there is an extension to a morphism $K^\triangleright \to \scI$~\cite[Def.~5.3.1.7]{Lurie:HTT}; that is, each $\kappa$-small diagram in $\scI$ admits a cocone.}
\cite[Def.~5.3.4.5]{Lurie:HTT}.
A morphism $Y \to X$ in $\scX$ is called \textit{relatively $\kappa$-compact}~\cite[Def.~6.1.6.4]{Lurie:HTT} if, for each morphism $A \to X$ from a $\kappa$-compact object $A$, the pullback $A \times_X Y$ is again a $\kappa$-compact object in $\scX$.
We let $S^\kappa$ denote the class of relatively $\kappa$-compact morphisms in $\scX$.
By the pasting law for pullbacks, the class $S^\kappa$ is closed under pullback.

We now want to say that there is a morphism $\Obj^\kappa_* \to \Obj^\kappa$ such that any relatively $\kappa$-compact morphism $Y \to X$ fits into a (unique) pullback square
\begin{equation}
\begin{tikzcd}
	Y \ar[d] \ar[r]
	& \Obj^\kappa_* \ar[d]
	\\
	X \ar[r]
	& \Obj^\kappa
\end{tikzcd}
\end{equation}
This is done as follows:
let $S$ be a class of morphisms in $\scX$ which is closed under pullback.
We let \smash{$\scO_\scX^S \subset \scFun(\Delta^1, \scX)$} denote the subcategory whose objects are the morphisms in $S$ and whose morphisms $(B \to A) \longrightarrow (Y \to X)$ are those commutative squares
\begin{equation}
\begin{tikzcd}
	B \ar[r] \ar[d]
	& Y \ar[d]
	\\
	A \ar[r]
	& X
\end{tikzcd}
\end{equation}
which are cartesian (the higher morphisms are as in $\scFun(\Delta^1, \scX)$)~\cite[Notation~6.1.3.4]{Lurie:HTT}.
The inclusion $\Delta^{\{1\}} \hookrightarrow \Delta^1$ induces a right fibration \smash{$\scO_\scX^S \to \scX$}.

\begin{definition}
\cite[Def.~6.1.6.1]{Lurie:HTT}
Let $S$ be a class of morphisms in $\scX$ which is closed under pullback.
A \textit{classifying morphism for $S$} is a final object of $\scO_\scX^S$.
\end{definition}

The existence of a classifying morphism for the class $S^\kappa$ of relatively $\kappa$-small morphisms is a consequence of the following theorem, which Lurie attributes to Rezk:

\begin{theorem}
\emph{\cite[Thm.~6.1.6.8]{Lurie:HTT}}
An $\infty$-category $\scX$ is an $\infty$-topos if and only if it has the following properties:
\begin{myenumerate}
\item $\scX$ is presentable,

\item colimits in $\scX$ are universal (see Definition~\ref{def:oo-topos}(2)), and

\item for each sufficiently large regular cardinal $\kappa$, there exists a classifying morphism $\Obj^\kappa_* \to \Obj^\kappa$ for the class $S^\kappa$ of relatively $\kappa$-compact morphisms in $\scX$.
\end{myenumerate}
\end{theorem}

Let $\kappa$ be a regular cardinal such that the classifying morphism $\Obj^\kappa_* \to \Obj^\kappa$ exists in $\scX$, and let $V \in \scX$ be an object such that the canonical morphism $V \to *$ is relatively $\kappa$-compact.
Then, there is a unique pullback square
\begin{equation}
\label{eq:class square for V --> *}
\begin{tikzcd}
	V \ar[r, "\overline{V^\dashv}"] \ar[d]
	& \Obj^\kappa_* \ar[d]
	\\
	* \ar[r, "V^\dashv"']
	& \Obj^\kappa
\end{tikzcd}
\end{equation}
The \v{C}ech nerve $\cC(V^\dashv)$ of the morphism $V^\dashv \colon * \to \Obj^\kappa$ classifying the object $V$ describes a loop object $\Omega_{V^\dashv} \Obj^\kappa$ together with its $\infty$-group structure (which any loop object in $\scX$ carries, see Example~\ref{eg:oo-group actions}(3)).
The following is a reformulation of~\cite[Def.~4.9]{NSS:oo-bundles_I}:

\begin{definition}
Let $\scX$ be an $\infty$-topos, let $\kappa$ be a sufficiently large regular cardinal, and let $V \to *$ be a relatively $\kappa$-compact morphism.
The \emph{automorphism $\infty$-group of $V$} is the group object
\begin{equation}
	\Aut(V) \coloneqq \cC (V^\dashv)
	\quad
	\in \Grp(\scX)\,.
\end{equation}
\end{definition}

Consequently, we obtain a classifying object
\begin{equation}
	\rmB \Aut(V) \coloneqq |\cC(V^\dashv)|
	\quad
	\in \scX
\end{equation}
for $\Aut(V)$-principal $\infty$-bundles in $\scX$.
Moreover, there is a canonical action of $\Aut(V)$ on $V$, encoded by the simplicial object
\begin{equation}
	V \dslash \Aut(V) \coloneqq \cC(\overline{V^\dashv})
	\quad
	\in \scFun(N\bbDelta^\opp, \scX)\,.
\end{equation}
Indeed, associated to any cartesian diagram
\begin{equation}
\begin{tikzcd}
	B \ar[d, "p"'] \ar[r]
	& Y \ar[d, "q"]
	\\
	A \ar[r]
	& X
\end{tikzcd}
\end{equation}
there is a cartesian diagram
\begin{equation}
\begin{tikzcd}[column sep={1.5cm,between origins}, row sep={1cm,between origins}]
	& A \ar[rr] \ar[dd, equal]
	& & X \ar[dd, equal]
	\\
	B \ar[rr, crossing over] \ar[dd] \ar[ur, "p"]
	& & Y \ar[ur, "q"']
	&
	\\
	& A \ar[rr]
	& & X
	\\
	A \ar[ur, equal] \ar[rr]
	& & X \ar[ur, equal]  \ar[from=uu, crossing over]
	&
\end{tikzcd}
\end{equation}
in \smash{$\scFun(N\bbDelta_{+, \leq 0}^\opp, \scX)$}, where we view the vertical edges as the objects in \smash{$\scFun(N\bbDelta_{+, \leq 0}^\opp, \scX)$}.
Forming the \v{C}ech nerve is the right Kan extension along the inclusion $\iota \colon N\bbDelta_{+, \leq 0}^\opp \hookrightarrow N\bbDelta_+$ (and then restricting along $N\bbDelta^\opp \hookrightarrow N\bbDelta_+^\opp$ to obtain a simplicial object).
Since the right Kan extension $\iota_*$ is a right adjoint, we obtain a cartesian diagram
\begin{equation}
\begin{tikzcd}
	\cC p \ar[r] \ar[d]
	& \cC q \ar[d]
	\\
	\sfc A \ar[r]
	& \sfc X
\end{tikzcd}
\end{equation}
in $\scFun(N\bbDelta^\opp, \scX)$, where in the bottom row we have the constant diagrams on $A$ and $X$, respectively.
Applied to the cartesian square~\eqref{eq:class square for V --> *}, this produces a cartesian square
\begin{equation}
\begin{tikzcd}
	V \dslash \Aut(V) \ar[r] \ar[d]
	& \Aut(V) \ar[d]
	\\
	\sfc V \ar[r]
	& \sfc *
\end{tikzcd}
\end{equation}
which establishes $V \dslash \Aut(V)$ as an $\Aut(V)$-action on $V$ in $\scX$ (Definition~\ref{def:group-action}).
One can now show:

\begin{proposition}
\label{st:all oo-bundles are associated}
\emph{\cite[Prop.~4.10]{NSS:oo-bundles_I}}
Let $\scX$ be an $\infty$-topos and $\kappa$ a regular cardinal such that the classifying morphism $\Obj^\kappa_* \to \Obj^\kappa$ exists.
Let $V \in \scX$ be an object such that $V \to *$ is relatively $\kappa$-compact.
Then, each $\infty$-bundle $p \colon E \to X$ in $\scX$ with fibre $V$ is associated to an $\Aut(V)$-principal $\infty$-bundle $P \to X$ in $\scX$ via the canonical action $V \dslash \Aut(V)$.
\end{proposition}

\begin{theorem}
\emph{\cite[Thm.~4.11]{NSS:oo-bundles_I}}
With $\scX$ and $V$ as in Proposition~\ref{st:all oo-bundles are associated}, there is a bijection
\begin{equation}
	\pi_0 \Bun_V(X) \simeq \rmH^1 \big( X; \Aut(V) \big)
	\simeq \pi_0 \scX\big( X, \rmB\Aut(V) \big)
	\simeq \pi_0 \Bun \big( X; \Aut(V) \big)\,.
\end{equation}
In particular, $\infty$-bundles on $X$ with fibre $V$ are classified by $\rmH^1 \big( X; \Aut(V) \big)$.
\end{theorem}


\section{Interaction of principal $\infty$-bundles with $\infty$-functors}
\label{sec:Geometric mps}


We now compare principal $\infty$-bundles in different $\infty$-topoi.
Let $\scX$ be an $\infty$-topos and $G$ a group object in $\scX$.
Recall from Definition~\ref{def:G-pr oo-Bun} that a $G$-principal $\infty$-bundle over an object $X \in \scX$ is a group action over $X$ such that the canonical morphism from the geometric realisation $|P \dslash G|$ to $X$ is an equivalence in $\scX$.

\begin{proposition}
\emph{\cite[Thm.~3.32]{Bunk:Pr-oo-Bundles}}
Let $\scX$ and $\scY$ be $\infty$-topoi and $f \colon \scX \to \scY$ an $\infty$-functor which preserves finite products and geometric realisations.
Let $X \in \scX$ and $G \in \Grp(\scX)$.
Then, $f$ maps group objects in $\scX$ to group objects in $\scY$, $G$-actions over $X$ in $\scX$ to $f(G)$-actions over $f(X)$ in $\scY$, and $G$-principal $\infty$-bundles over $X$ in $\scX$ to $f(G)$-principal $\infty$-bundles over $f(X)$ in $\scY$.
\end{proposition}

One encounters such $\infty$-functors, in particular, in the case of \textit{cohesive} $\infty$-topoi.
This notion was introduced by Schreiber~\cite[Def.~4.1.8]{Schreiber:DCCT} as an $\infty$-categorical enhancement of ideas by Lawvere~\cite{Lawvere:Axiomatic_cohesion}.
We now recall this notion.
First, there is a type of $\infty$-functor between $\infty$-topoi which is particularly well-adapted to the additional structure present in $\infty$-topoi:

\begin{definition}
\cite[Def.~6.3.1.1]{Lurie:HTT}
Let $\scX$ and $\scY$ be $\infty$-topoi.
A \emph{geometric morphism $\scX \to \scY$} is an $\infty$-functor $f_* \colon \scX \to \scY$ which has a left-exact%
\footnote{An $\infty$-functor between finitely complete $\infty$-categories is left-exact if it preserves finite limits~\cite[Rmk.~5.3.2.3]{Lurie:HTT}.}
left adjoint $f^* \colon \scY \to \scX$.
\end{definition}

One can show that the $\infty$-topos $\scS$ is a final object in the $\infty$-category of $\infty$-topoi and geometric morphisms~\cite[Prop.~6.3.4.1]{Lurie:HTT}.
Thus, for each $\infty$-topos $\scX$, there is a canonical adjoint pair
\begin{equation}
\begin{tikzcd}
	\delta : \scS \ar[r, shift left=0.15cm, "\perp"' yshift=0.05cm]
	& \scX : \Gamma\,, \ar[l, shift left=0.15cm]
\end{tikzcd}
\end{equation}
whose right adjoint $\Gamma$ is a geometric morphism.
It is often called the \textit{global sections $\infty$-functor} of $\scX$.

\begin{definition}
\cite[Def.~4.1.8]{Schreiber:DCCT}
An $\infty$-topos is \textit{cohesive} if its global sections $\infty$-functor is part of a triple adjunction $\Pi \dashv \delta \dashv \Gamma \dashv \mathrm{codisc}$, satisfying that $\delta$ and $\mathrm{codisc}$ are fully faithful and $\Pi$ preserves finite products.
\end{definition}

Important examples of cohesive $\infty$-topoi consist of the $\infty$-topoi $\scH = \PSh(N\Cart)$ and $\scH_\tau = \Sh(N\Cart, \tau)$ from Example~\ref{eg:oo-topoi of smooth spaces} (see~\cite[Prop.~4.1.32]{Schreiber:DCCT}).
In particular, under the presentations in~\eqref{eq:presentation of PSh} and~\eqref{eq:presentation of Sh} the global-section adjunction arises from the adjunction $const \dashv \ev_{\RN^0}$~\cite[Prop.~4.1.30]{Schreiber:DCCT}.
It follows that there is a canonical equivalence $\Pi \simeq \colim$ for the additional left adjoint $\infty$-functor.
If one interprets objects of $\scH$ or $\scH_\tau$ as higher smooth spaces, $\Pi$ has an interpretation as taking an \textit{underlying space}, or a smooth version of the \textit{singular complex} functor in topology (see~\cite{BEBdBP:Class_sp_of_oo-sheaves, Bunk:Pr-oo-Bundles, ADH:Diff_Coho, Bunk:Sm_Spaces, Pavlov:DfgSp_and_Oka_Principle} for more background).

\begin{corollary}
The $\infty$-functor $\Pi \colon \scH_\tau \to \scS$ preserves group actions and principal $\infty$-bundles.
\end{corollary}


\section{Outlook: higher connections and characteristic classes}
\label{sec:Connections}


Once a good notion of $\infty$-bundles is in place, a crucial next step is to establish a theory of $\infty$-connections on these bundles.
This vast and important theory is currently still incomplete and under active development, and we only touch upon some of the current directions in this section.

For particularly well-understood examples of principal $\infty$-bundles, the 2-bundles (i.e.~where the structure group is 2-truncated), a full notion of connections exist; see, for instance,~\cite{Waldorf:Conns_on_pr_2-Bdls, Waldorf:PT_in_pr_2-Bdls}, going back to~\cite{BS:HGT, BH:Invitation_to_HGT} (but see also~\cite{FMP:2D_holonomy, Kapranov:Membranes} for approaches to two- and higher-dimensional parallel transport).
However, even in these cases it appears that one can emphasise different aspects of ordinary connections from differential geometry in the process of categorification.
Focussing on the parallel transport aspect leads to connections satisfying the \textit{fake curvature condition}, which stems from the interchange law in the path 2-groupoid of a manifold~\cite{BS:HGT}.
However, this condition appears too strong in certain situations:
for instance, if a manifold $M$ admits a String structure~\cite{Killingback:WS_Anomalies, ST:What_is_an_elliptic_object, Waldorf:String_and_Chern-Simons} (and see also Section~\ref{sec:Appls in physics} below), it admits a connection satisfying the fake curvature condition only if the tangent bundle $TM$ admits a \textit{flat} connection.
There has recently been a proposal to alleviate this by an \textit{adjustment} to the notion of a 2-connection~\cite{SS:M5-brane_models_II}.

For other simple structure groups, such as iterated deloopings of abelian Lie groups in $\scH_\tau$, there exists a full theory of $\infty$-connections via the Dold-Kan correspondence and Deligne complexes; see Example~\ref{eg:pr-oo-bdls from classifying objects}.
These examples of connections are particularly well understood.
For instance, for connections on $\infty$-bundles in $\scH_\tau$ classified by the objects $\rmB^{n-k} \rmB_\nabla^k \rmU(1)$, for $0 \leq k < n \in \NN$, a theory of moduli $\infty$-stacks for solutions to higher-gauge theoretic equations has recently been developed in~\cite{BS:Higher_Syms_and_Deligne_Coho}.

In a different direction, $\infty$-connections have been studied on principal $\infty$-bundles whose structure group arises as an integration of an $L_\infty$-algebra $\frg$~\cite{FSS:Cech_diff_char_classes_via_L_infty, Schreiber:DCCT}.
In particular, that facilitated the construction of (differential) characteristic classes and Chern-Weil theory~\cite{FSS:Cech_diff_char_classes_via_L_infty} in these situations (see also~\cite{FH:CW_forms_and_abstract_HoThy} for the case of Lie groups).
This has recently been developed further in close relation with rational homotopy theory in~\cite{FSS:Character_map}:
by Section~\ref{sec:Non-ab coho and BG} each connected object $T$ in an $\infty$-topos $\scX$ induces a non-abelian cohomology theory via
\begin{equation}
	\rmH^0(X;T) = \pi_0 \scX(X,T)
	\simeq \pi_0 \scX(X, \rmB \Omega T)
	\simeq \pi_0 \Bun(X; \Omega T)\,.
\end{equation}
For $\scX = \scS$, rational homotopy theory associates to each connected, nilpotent, rationally finite space $T \in \scS$ an $L_\infty$-algebra $\frl T$.
This controls the characteristic classes in real cohomology associated to $\Omega T$-principal $\infty$-bundles, or, equivalently, non-abelian cohomology with coefficients in $T$~\cite[Ch.~4]{FSS:Character_map}.
In particular, the theory of $\infty$-connections provides a starting point for the development of non-abelian differential cohomology theories and its geometric cocycles.


\section{Some applications in physics}
\label{sec:Appls in physics}


We conclude with a necessarily incomplete selection of occurrences of higher principal bundles in mathematical physics (in addition to those touched upon in Example~\ref{eg:pr-oo-bdls from classifying objects}).

\paragraph*{The B-field in string theory.}

One of the first examples of connections on higher bundles in mathematical physics appeared in string theory.
More concretely, it was realised that the B-field is captured mathematically by a connection on a \textit{gerbe}~\cite{Kapustin:D-branes_in_nontriv_B-fields} (see also~\cite{Murray:Bundle_gerbs, Murray:Intro_to_gerbes, Brylinski:Loops_and_GeoQuan, Bunk:Gerbes_in_Geo_and_FT} for background on gerbes).
As pointed out in Example~\ref{eg:pr-oo-bdls from classifying objects}, a gerbe can also be described as a particular model for a principal $\infty$-bundle in $\scH_\tau$ with structure group $\rmB \rmU(1)$.
This perspective was developed in~\cite{FSS:Cech_diff_char_classes_via_L_infty, Schreiber:DCCT}, building on earlier work~\cite{Gajer:Geo_of_Deligne_Coho}.

The statement that the B-field term in the string world-sheet actions---the Wess-Zumino-Witten term---and its associated Chan-Paton terms are modelled by connections on gerbes and vector bundles twisted by gerbes has been made precise in~\cite{BW:Transgression_of_D-branes, BW:OCTFTs_and_Gerbes} by enhancing these world-sheet actions, and thus the two-dimensional holonomy of gerbe connections and Chan-Paton bundles on D-branes, into a smooth functorial field theory in the sense of Stolz-Teichner~\cite{ST:SuSy_FTs_and_generalised_coho}.

\paragraph*{String groups and string structures.}

Each compact, simple and simply connected Lie group $G$ satisfies $\pi_2(G) = 0$, $\pi_3(G) \cong \ZN$, and thus $\rmH^3(G;\ZN) \cong \ZN$.
In~\cite{Stolz:Pos_Ric} Stolz proved, by giving an explicit construction, that there exists a morphism of topological groups $p \colon \String(G) \to G$ which, as a continuous map, is a 3-connected cover of $G$.
That is, $\pi_3(\String(G)) \cong 0$ and $\pi_i p$ is an isomorphism for all $i \neq 3$.
The codomain of any such map is called a \textit{string group extension of $G$}.
It is impossible to construct $\String(G)$ as a finite-dimensional manifold (the fibre of $p$ has cohomology in each even degree).
In particular, one writes $\String(n)$ for the case where $G = \Spin(n)$ is a spin group with $n \geq 3$.
A \textit{string structure} on a spin manifold $M$ is a lift of its $\Spin(n)$-principal bundle to a $\String(G)$-principal bundle, possibly in a higher geometric sense.

String structures, and thus string groups, feature in the differential geometry of free loop spaces $LM$ of manifolds $M$:
spin structures on $LM$ are equivalent to string structures on $M$~\cite{Killingback:WS_Anomalies, Waldorf:String_v_Spin}.
It has been a long-standing open problem to make sense of the hypothetical Dirac operator on $LM$, going back to~\cite{Witten:Elliptic_general, Witten:Loop_index} (see also~\cite[Sec.~2]{Stolz:Pos_Ric} for a review).
In order to carry out such differential geometric constructions, it is necessary to have a smooth geometric enhancement of the topological group $\String(G)$.
This has been achieved in the settings of $\infty$-dimensional Lie groups~\cite{NSW:Smooth_String_Model}, as well as 2-groups in various $\infty$-categories of smooth spaces (see, for instance,~\cite{SP:String, BCSS:String_via_Loops, BMS:Sm2Grp, Bunk:Pr-oo-Bundles, FRS:U1-Gerbe_connections, Waldorf:String_via_TR, Waldorf:Multiplicative_Gerbes}).

Via $\infty$-bundles, one can formulate lifts of (ordinary) principal $G$-bundles on $M$ to $\String(G)$-bundles as follows.
We work in the $\infty$-topos $\scH_\tau$ of $\infty$-sheaves on cartesian spaces.
There is a canonical generator of $\rmH^3(G;\ZN)$.
Geometrically, it is presented by the \textit{basic gerbe} on $G$~\cite{Meinrenken:Basic_Gerbe}.
This, in turn, is a $\rmB\rmU(1)$-principal $\infty$-bundle on $G$, classified by a morphism%
\footnote{In other words, Meinrenken's construction of the basic gerbe proves that the morphism $G \to \rmB^2\rmU(1) \simeq K(\ZN,2)$ in $\scS$ admits an enhancement to a morphism in $\scH_\tau$.}
$c_\CG \colon G \to \rmB^2\rmU(1)$ in $\scH_\tau$.
Waldorf showed~\cite{Waldorf:Multiplicative_Gerbes} that the basic gerbe admits a lift of the group structure of $G$, and the resulting group object is a string group extension for $G$.
We obtain a pullback square
\begin{equation}
\begin{tikzcd}
	\String(G) \ar[r] \ar[d]
	& * \ar[d]
	\\
	G \ar[r, "c_\CG"']
	& \rmB^2 \rmU(1)
\end{tikzcd}
\end{equation}
in $\scH_\tau$ which exhibits an extension of $\infty$-groups (see also Definition~\ref{def:oo-group extension}).
Thus, we also have a fibre sequence of classifying objects in $\scH_\tau$,
\begin{equation}
\begin{tikzcd}
	\rmB\String(G) \ar[r] \ar[d]
	& * \ar[d]
	\\
	\rmB G \ar[r, "\rmB c_\CG"']
	& \rmB^3 \rmU(1)
\end{tikzcd}
\end{equation}
Consider an ordinary principal $G$-bundle $P \to M$, classified by a morphism $c_P \colon M \to \rmB G$ in $\scH_\tau$.
The space of string structures for this principal $G$-bundle is the space---in fact, the groupoid---of lifts of $c_P$ through $\rmB\String(G)$.
By the universal property of pullbacks in the $\infty$-topos $\scH_\tau$, this groupoid is equivalent to that of trivialisations of the composition $(\rmB c_\CG) \circ c_P \colon M \to \rmB^3\rmU(1)$.
This composite classifies a $\rmB^2\rmU(1)$-principal $\infty$-bundle on $M$, which is also known as the \textit{Chern-Simons 2-gerbe} of the principal $G$-bundle $P \to M$~\cite{CJMSW:BGrbs_for_CS_and_WZW}.
In other words, from the $\infty$-bundle perspective, string structures for a principal $G$-bundle $P \to M$ are equivalent to trivialisations of the Chern-Simons 2-gerbe of $P \to M$~\cite[Def.~1.1.5]{Waldorf:String_and_Chern-Simons}.

\paragraph*{Cohomotopy and hypothesis H.}

It is a widely accepted paradigm that charges of D-branes in string theory are classified by (twisted) K-theory~\cite{Witten:D-branes_and_KT, MM:KT_and_RR_charge}.
A recent proposal for a similar framework for charge quantisation in M-theory is known as \textit{Hypothesis H} and goes back to~\cite{Sati:Framed_M-branes_and_top_invars}; we refer the reader to~\cite{FSS:Rational_HStr_of_M-theory} for a review and further references.
The proposal rests on the observation that the differential form data of M-theory matches, via the character map mentioned in Section~\ref{sec:Connections}, with the characters obtained from (twisted) non-abelian cohomology with coefficients in the 4-sphere $\bbS^4$ (in the sense of Section~\ref{sec:Non-ab coho and BG}).
That is, at the level of spaces it is controlled by the $\infty$-functors
\begin{equation}
	\scS(-, \bbS^4) \simeq \scS(-, \rmB \Omega \bbS^4)
	\simeq \Bun(-; \Omega \bbS^4)\,.
\end{equation}
It can thus be described equivalently as the study of cohomotopy theory with coefficients in $\bbS^4$, or as the study of principal $\infty$-bundles for the based loop group of $\bbS^4$ in the sense of Section~\ref{sec:Motivation} (see also~\cite[Ex.~2.10]{FSS:Character_map}).

\begin{small}

\makeatletter

\interlinepenalty=10000

\makeatother

\bibliographystyle{alphaurl}
\addcontentsline{toc}{section}{References}
\bibliography{Infty-Bundles_Bib}

\vspace{0.5cm}

\noindent
Mathematical Institute, The University of Oxford.
\\
severin.bunk@maths.ox.ac.uk

\end{small}

\end{document}

%% file: Infty-Bundles_Defs.tex
\newcommand{\Obj}{{\mathrm{Obj}}}
\newcommand{\im}{{\mathrm{im}}}
\newcommand{\pr}{\mathrm{pr}}

\newcommand{\Fun}{{\mathrm{Fun}}}

\newcommand{\Spin}{{\mathrm{Spin}}}
\newcommand{\String}{{\mathrm{String}}}

\newcommand{\Sp}{{\mathrm{Sp}}}

\newcommand{\rmH}{\mathrm{H}}
\newcommand{\rmU}{{\mathrm{U}}}

\newcommand{\rmB}{{\mathrm{B}}}
\newcommand{\rmE}{{\mathrm{E}}}

\newcommand{\ev}{\operatorname{ev}}

\newcommand{\dd}{\mathrm{d}}

\newcommand{\id}{{\mathrm{id}}}

\newcommand{\opp}{\mathrm{op}}

\newcommand{\colim}{{\mathrm{colim}}}

\newcommand{\Aut}{\mathrm{Aut}}
\renewcommand{\lim}{{\mathrm{lim}}}

\newcommand{\scO}{{\mathscr{O}}}
\newcommand{\scX}{{\mathscr{X}}}
\newcommand{\scI}{\mathscr{I}}
\newcommand{\scK}{{\mathscr{K}}}
\newcommand{\scH}{\mathscr{H}}

\newcommand{\scC}{\mathscr{C}}
\newcommand{\scD}{\mathscr{D}}

\newcommand{\scY}{{\mathscr{Y}}}

\newcommand{\scS}{\mathscr{S}}

\newcommand{\CG}{\pazocal{G}}
\newcommand{\CU}{{\pazocal{U}}}

\newcommand{\sfc}{\mathsf{c}}

\newcommand{\NN}{\mathbb{N}}
\newcommand{\RN}{\mathbb{R}}
\newcommand{\ZN}{\mathbb{Z}}

\newcommand{\bbS}{\mathbb{S}}

\newcommand{\bbE}{\mathbb{E}}

\newcommand{\frg}{{\mathfrak{g}}}
\newcommand{\frl}{\mathfrak{l}}

\newcommand{\cC}{{\check{C}}}

\newcommand{\Grb}{{\mathscr{G}\mathrm{rb}}}

\newcommand{\Set}{{\mathscr{S}\mathrm{et}}}
\newcommand{\sSet}{{\mathscr{S}\mathrm{et}_{\hspace{-.03cm}\Delta}}}
\newcommand{\scFun}{{\mathscr{F}\hspace{-.03cm}\mathrm{un}}}

\newcommand{\Cat}{{\mathscr{C}\mathrm{at}}}

\newcommand{\Mfd}{{\mathscr{M}\mathrm{fd}}}

\newcommand{\Cart}{{\mathscr{C}\mathrm{art}}}

\newcommand{\Sh}{{\mathscr{S}\mathrm{h}}}
\newcommand{\PSh}{{\mathscr{PS}\mathrm{h}}}
\newcommand{\Top}{{\mathscr{T}\mathrm{op}}}
\newcommand{\Gpd}{{\mathscr{G}\mathrm{pd}}}
\newcommand{\Grp}{{\mathscr{G}\hspace{-0.02cm}\mathrm{rp}}}

\newcommand{\Bun}{{\mathscr{B}\hspace{-0.02cm}\mathrm{un}}}

\newcommand{\eq}{\overset{\simeq}{\longrightarrow}}

\newenvironment{myenumerate}{\begin{enumerate}[topsep=-\parskip+0.1cm, itemsep=0.1cm, parsep=0cm, leftmargin=*, label=(\arabic*)]}{\end{enumerate}}

\newcommand{\qen}{\hfill$\triangleleft$}

\newcommand{\dslash}{{/\hspace{-0.1cm}/}}

\newcommand{\bbDelta}{{\mathbbe{\Delta}}}